%% file: Dis.vectorielles2001.tex
\overfullrule=0mm
\def\det{{\rm det}}
\def\d¡{{\rm d¡}}
\def\tr{{\rm tr}}

\def\min{{\rm min}}
\def\sup{{\rm sup}}
\def\adots{\mathinner{\mkern2mu\raise1pt\hbox{.}
\mkern3mu\raise4pt\hbox{.}\mkern1mu\raise7pt\hbox{.}}}

\def\Pf{{\rm Pf}}
\def\Log{{\rm Log}}

\def\GL{{\rm GL }}
\def\Det{{\rm Det}}
\input AMSsym.def
\input AMSsym.tex

\input Init.Som
\def\L{{\rm L }}
\def\P{{\rm P }}
\def\PSL{{\rm PSL }}
\def\SL{{\rm SL }}

\def\Pf{{\rm Pf}}
\def\SO{{\rm SO}}

\catcode`\@=12

 \centerline{\bf Distributions vectorielles homog\`enes sur une alg\`ebre de Jordan}
 
\bigskip

\bigskip

\bigskip

\centerline{\bf Bruno Blind}

\centerline{\it Institut Elie Cartan, U.M.R 7502}

\centerline{\it Universit\'e H.Poincar\'e, B.P.239, 54506 Vandoeuvre-l\'es-Nancy.}

\bigskip

\bigskip

\rightline{\it{Publish or perish.}}

\bigskip
\rightline{\it{A Marraine.}}

\rightline{\it{"On reste l\`a, il fait si bon parfois quand le soir tombe}}

\rightline{\it{et qu'on regarde simplement ses mains."}}

\rightline{C. Esteban}

\bigskip

\bigskip  

\bigskip

\bigskip

\footnote{}{e-mail: blind@iecn.u-nancy.fr}

\footnote{}{tel: 0383684587; fax: 0383684534}

\bigskip

\vfill \eject 

\bigskip

{\bf Abstract:} We study distributions on a Euclidean Jordan algebra $V$ with values in a finite
dimensional representation space for the identity component $G$ of the structure group of $V$ and
homogeneous equivariance condition. We show that such distributions exist if and only if the
representation is spherical, and that then the dimension of the space of these distributions is
$r+1$ ( where $r$ is the rank of $V$). We give also construction of these distributions and of those that are
invariant under the semi-simple part of $G$.

\bigskip

\bigskip  

\bigskip

\bigskip

\footnote{}{Classification AMS: 46F10,17C.}

\footnote{}{Mots clefs: Distributions homog\`enes, Distributions vectorielles, Alg\`ebres de
  Jordan.}

\vfill \eject

\beginsection {I. Introduction.}

\bigskip

Soient $V$ une alg\`ebre de Jordan r\'eelle, simple
et euclidienne et $G$ la composante neutre de son groupe de structure (voir le paragraphe II pour les
rappels sur les alg\`ebres de Jordan; on peut ici prendre pour $V$ l'espace
$\hbox{Sym}(r, {\Bbb R})$ des matrices sym\'etriques r\'eelles d'ordre
$r$, et pour $G$ le groupe $\SL (r,{\Bbb R})$
agissant par: $g.X = gX{}^tg$). Soit d'autre part
$(\pi, L_{\pi})$ une repr\'esentation irr\'eductible
de dimension finie de $G$. Le but de cet article est l'\'etude 
 des distributions homog\`enes sous
l'action de $G$ \`a valeurs dans l'espace  $L_{\pi}$. Nous montrons, dans les
paragraphes III et V, que de telles distributions non nulles existent si et seulement si la
repr\'esentation
$\pi$ est sph\'erique, et que dans ce cas la dimension de l'espace des distributions
homog\`enes de degr\'e donn\'e est
\'egal aux nombres d'orbites ouvertes du groupe
$G$, r\'esultat bien connu dans le cas scalaire (cf [Ri.-St.], [Mu.] par exemple). Nous
construisons, au paragraphe IV, ces distributions homog\`enes. Nous d\'eterminons aussi, au
paragraphe V, l'espace des distributions
$G^1$-invariantes (o\`u
$G^1$ est la "partie semi-simple" du groupe $G$) port\'ees par le lieu singulier de
l'alg\`ebre de Jordan. 

Nos d\'emonstrations s'inspirent de celles \'elabor\'ees par F. Ricci et E.M.Stein dans
leur travail ([Ri.-St.]) consacr\'e aux distributions homog\`enes sur
les matrices hermitiennes; en particulier la d\'emonstration de la caract\'erisation des
repr\'esentations, obtenue au paragraphe III, est une adaptation du lemme 5 de  [Ri.-St.] (dans
[Bl.1], nous avions d\'eja
\'etendu ce lemme au cas scalaire pour toute alg\`ebre de Jordan simple et euclidienne). 

Ce travail constitue une g\'en\'eralisation naturelle de certains
r\'esultats que J.A.C.Kolk et V.S. Varadarajan ont obtenus 
pour les distributions $\SO_o(1,n)$-invariantes sur ${\Bbb R}^{1+n}$ \`a valeurs dans l'espace d'une
repr\'esentation irr\'eductible de dimension finie du groupe $\SO_o(1,n)$, et 
\`a support dans la nappe sup\'erieure du c\^one de lumi\`ere (voir [Kol.-Va.1] ainsi que
[Kol.-Va.2]).

Signalons aussi l'article de J.Hilgert et K.H.Neeb ([Hi.-Ne.]), o\`u les auteurs \'etudient les
distributions de Riesz vectorielles sur les alg\`ebres de Jordan euclidiennes.

Je remercie J.L.Clerc de m'avoir signal\'e le travail de J.A.C.Kolk et V.S. Varadarajan, H.
Rubenthaler de ses (nombreux) encouragements amicaux ainsi que P.Y. Gaillard et F. Chargois pour
les discussions que nous avons eues ensemble.

\vfill \eject

\beginsection {II. Rappels sur les alg\`ebres de Jordan.}

\bigskip

\noindent $\underline {\hbox {A. Les alg\`ebres de Jordan simples et euclidiennes.}}$

\bigskip

Nous renvoyons \`a [Br.-Koe.], [Fa.-Ko.], [Sa.1] pour les d\'efinitions et les principaux faits
concernant les alg\`ebres de Jordan, et reprenons ici les rappels de [Bl].

Soit $V$  une alg\`ebre de Jordan euclidienne sur le corps  
${\Bbb R}$, simple, de dimension n, de rang r et
d'\'el\'ement unit\'e e. Pour $x$ et $y$ dans $V$, on d\'efinit
les endomorphismes de $V$:

\noindent $\L (x)y = xy$, $\P(x) = 2\L^2(x)-\L(x^2)$ et $x \square y = \L(xy) + [\L(x),\L(y)]$.

\noindent On munit $V$ de la forme bilin\'eaire d\'efinie
positive $< x,y > = {r\over n}\tr \L(xy)$; si $g$ est un \'el\'ement de $\GL(V)$, on note par
$g^*$ son adjoint par rapport \`a ce produit scalaire; par ailleurs $dx$ d\'esignera la mesure
euclidienne associ\'ee au produit $<.,.>$.

On note $\det x $ le d\'eterminant de $V$. C'est un polyn\^ome irr\'eductible
et homog\`ene de degr\'e $r$. Un \'el\'ement $x$ de $V$ est inversible si $\det x \not=   0$. Nous
noterons $\Det g$ le d\'eterminant d' un \'el\'ement $g$ de $\GL (V)$, et $V^{\times}$ l'ensemble
des \'el\'ements inversibles de l'alg\`ebre $V$.

 Tout \'el\'ement $x$ de l'alg\`ebre $V$ peut s'\'ecrire sous la forme: $x = \sum_{i
=0}^r{\lambda_ie_i},  \lambda_i \in {\Bbb R}$ o\`u $\{e_1,\ldots ,e_r\}$ est un syst\`eme complet
d'idempotents primitifs orthogonaux (d\'ependant de $x$). Les nombres $\lambda_i$ sont appel\'es les
valeurs propres de $x$ et l'on a: $\det x = \prod_{i=0}^r \lambda_i$. Le rang de $x$ est le nombre de
r\'eels $\lambda _i$ non nuls (compt\'es avec multiplicit\'e). Si un \'el\'ement inversible $x$
poss\`ede $p$ valeurs propres positives et $q$ valeurs propres n\'egatives, on dira que sa signature
est $(p,q)$ ($p+q = r$). On note $\Omega_j$ l'ensemble des \'el\'ements inversibles de signature
$(r-j,j)$. Le c\^one $ \Omega = \Omega_0$ est la composante connexe contenant e de
l'ensemble des \'el\'ements inversibles de $V$.

Notons $G$ la composante neutre du groupe $G(\Omega)$ des isomorphismes lin\'eaires de $V$
pr\'eservant le c\^one $\Omega$ (c'est aussi la composante neutre du groupe de structure de $V$); on
a la d\'ecomposition
$G = {\Bbb R}^{+*}\times G^1$ o\`u
$G^1 = \{ g\in G / \Det g = 1 \}$ est un groupe de Lie connexe, semi-simple de centre trivial.
Fixons un syst\`eme complet d'idempotents orthogonaux $\{e_1,\ldots ,e_r\}$ de
l'alg\`ebre $V$. Le groupe $G$ 
admet une d\'ecomposition d'Iwasawa $G = KAN$ o\`u $K$ est le
sous-groupe de
$G$ des \'el\'ements qui fixent l'unit\'e e, $A$ est engendr\'e par les $\P (\sum_{i=1}^r a_ie_i)$
(avec $a_i >0$), et o\`u $N$ est engendr\'e par les transformations de Frobenius $\exp(2z
\square e_j)$ avec $z$ dans $ V_{ij} = \{x\in V /e_ix = e_jx = {1\over 2}x\}$,$i<j$ (pour tout ceci
voir [Fa.-Ko.]); rappelons que tous les espaces $ V_{ij}$ sont de m\^eme dimension $d$ et que
l'on a $n = r + {dr(r-1)\over2} $. Notons que ce choix de
$N$ n'est pas celui de [Fa.-Ko.]: notre $ N$ correspond \`a leur $\bar {N}$.

Les orbites de $V$ sous l'action du groupe
$G$ ont
\'et\'e d\'etermin\'ees par S.Kaneyuki et I.Satake (cf.[Ka.], [Sa.2]): ce sont les
ensembles
$S_{p,q} = Go_{p,q}$ o\`u $o_{p,q} = \sum_{i=1}^{p-q}{e_i} -
\sum_{i=1}^q{e_{p-q+i}}$ (avec $0 \leq q \leq p \leq r$); en particulier la signature est invariante
sous
$G$. Pour
$p = r$ les orbites
$S_{p,q}$ sont ouvertes et l'on a : $S_{r,q}$ = $\Omega_q$; les
autres orbites sont contenues dans le lieu singulier $S = \{ x
\in V \ / \ \det x = 0 \}$  de l'alg\`ebre $V$ et sont en fait des $G^1$
orbites. Pour une telle orbite singuli\`ere, notons $T_{pq}$ (resp. $N_{pq}$)
l'espace tangent (resp. l'espace normal) en $o_{p,q}$. On v\'erifie alors facilement que si 
$$V = \bigoplus _{1\leq i\leq j\leq r}V_{ij}$$ est la d\'ecomposition de Peirce associ\'ee au
syst\`eme $\{e_1,\ldots ,e_r\}$,
l'espace $N_{pq}$ s'identifie \`a:
$\bigoplus _{p < i\leq j\leq r}V_{ij}$
et $T_{pq}$ \`a:
$\bigoplus _{1\leq i\leq j\leq r \atop i \leq p}V_{ij}$.

\noindent En particulier, $N_{pq}$ est une alg\`ebre de Jordan simple
et euclidienne de rang
$r-p$ que nous noterons $V_{(r-p)}$.
Remarquons que $W_{pq} = \bigcup_{p'=p}^r\bigcup_{q'=q}^{q+p'-p}S_{p'q'}$, est un voisinage
ouvert de l'orbite $S_{pq}$, invariant sous $G$.

Soit $V^{\Bbb C}$ la complexifi\'ee de l'alg\`ebre $V$, le groupe $G_{\Bbb C}$, composante neutre
du groupe de structure de $V^{\Bbb C}$, est un complexifi\'e de $G$, au sens o\`u
son alg\`ebre de Lie est une complexifi\'ee de l'alg\`ebre de Lie  ${\goth G}$ de $G$. On note
$G^1_{\Bbb C}$ le sous-groupe connexe de
$G_{\Bbb C}$ d'alg\`ebre de Lie ${\goth G^1}_{\Bbb C}$ complexifi\'ee de l'alg\`ebre de Lie ${\goth
G^1}$ de $G^1$.

 Rappelons bri\`evement la classification des alg\`ebres de Jordan 
sur
${\Bbb R}$ simples et euclidiennes ainsi que de leur complexifi\'e([Br.Koe.], [Fa.-Ko.], [Sp.]), on
donne aussi l'indice $(G(\Omega) : G)$ de $G$ dans $G(\Omega)$:

1) $V = {\Bbb R}$, le groupe $G^1$ est ici trivial.

2) $V = \hbox{Sym}(r, {\Bbb R})$( $r > 1$), l'espace des matrices sym\'etriques r\'eelles d'ordre
$r$ avec
$$x.y = {1\over 2}(xy +yx) \qquad d = 1 \qquad n = {r(r+1)\over 2} \qquad < x,y > = \tr xy$$
Ici le d\'eterminant est le d\'eterminant usuel. Le groupe $G^1$ est $\PSL (r,{\Bbb R})$ agissant par: $g.X =
gX{}^tg$. On a $V^{\Bbb C} = \hbox{Sym}(r, {\Bbb C})$ et $G^1_{\Bbb C}$ est le groupe des transformations
lin\'eaires $g.X = gX{}^tg$ avec $g$ dans $\SL (r,{\Bbb C})$, c'est \`a dire s'identifie \`a $SL (r,{\Bbb C})/
\{Id,(-1)^{r+1}Id \}$. On a
$(G(\Omega) : G) = 2$ si $n$ est pair et $1$ sinon.

3) $V = \hbox{Herm}(r,{\Bbb C})$( $r > 1$), l'espace des matrices hermitiennes complexes d'ordre
$r$ avec
$$x.y = {1\over 2}(xy +yx) \qquad d = 2 \qquad n = r^2 \qquad < x,y > = \tr xy$$
Ici encore le d\'eterminant est le d\'eterminant usuel. Le groupe $G^1$ est $\PSL (r,{\Bbb C})$
agissant par: $g.X = gX{}^t\bar g$. On a $V^{\Bbb C} = \hbox{Mat}(r, {\Bbb C})$ et $G^1_{\Bbb C}$
est le groupe des transformations lin\'eaires $(g_{1},g_{2}).X = g_{1}Xg_{2}^{-1}$ avec $g_{1}$ et
$g_{2}$ dans
$\SL (r,{\Bbb C})$. Enfin $(G(\Omega) : G) = 2$.

4) $V = \hbox{Herm}(r,{\Bbb H})$( $r > 2$), l'espace des matrices hermitiennes sur $\Bbb H$
d'ordre $r$, avec le produit $x.y = {1\over 2}(xy +yx)$. Ici $d = 4$. On peut aussi voir
l'alg\`ebre $V$ dans l'espace des matrices antisym\'etriques $\hbox{Alt}(2r, {\Bbb C})$: $V = \{x
\in \hbox{Alt}(2r, {\Bbb C}) /  J\bar x = xJ \}$ avec la multiplication $x.y = {1\over 2}(xJy
+yJx)$ o\`u $J$ est la matrice $$J = \pmatrix{ 0 & I_r \cr
-I_r & 0 \cr
}.$$
Le groupe $G^1$ est $\PSL (r,{\Bbb H})$
agissant par: $g.X = gX{}^tg$.

 On a $V^{\Bbb C} =
\hbox{Alt}(2r, {\Bbb C})$ et
$G^1_{\Bbb C}$ est le
groupe des transformations lin\'eaires $g.X = gX{}^tg$ avec $g$ dans $\SL (2r,{\Bbb C})$. Le
d\'eterminant est donn\'e par le Pfaffien: $\det x = \Pf (JxJ)$ et $(G(\Omega) : G) = 1$.

5) $V = {\Bbb R}\times E $, o\`u $E$ est un espace euclidien de produit scalaire $(.,.)_E$,
de dimension $n-1$,avec $n > 4$. Le produit de Jordan est d\'efini par: $(\lambda,u)(\mu,v) =
(\lambda
\mu + (u,v)_E ,\lambda v + \mu u)$; ici $r = 2$, $d = n-2$ et le d\'eterminant est donn\'e par:
$\det (\lambda,u) =
\lambda^2 - \Vert u \Vert^2 _E$.Le groupe $G^1$ est la composante neutre de $\SO (1,n-1)$ agissant
naturellement sur $V$. On a $V^{\Bbb C} = {\Bbb C}\times {{\Bbb C}}^{n-1}$ et
$G^1_{\Bbb C} = \SO (n,{\Bbb C})$. $(G(\Omega) : G) = 2$.

6) L'alg\`ebre de Jordan exceptionnelle $\hbox{Herm}(3,{\Bbb O})$ constitu\'ee de matrices $(3,3)$
\`a coefficients dans les octaves de Cayley ($d = 8$, $r = 3$, $n = 27$), dans ce cas on sait que
le groupe $G^1_{\Bbb C}$ est le groupe complexe simplement connexe $E_6$. $(G(\Omega) : G) = 1$.

\bigskip 

Quelques remarques pour terminer ce paragraphe. Si $x$ est un \'el\'ement inversible de $V$ tel que
$\P (x)e = e$, $\P (x)$ appartiendra \`a $G(\Omega)$; c'est en particulier le cas des $g_i = \P(e_1 +
\ldots - e_i+ \ldots +e_r)$ pour $1\leq i \leq r$. Il existe un automorphisme $w_i$ de l'alg\`ebre
$V$ tel que $w_i(-e_1 + \ldots +e_r) = e_1 + \ldots - e_i+ \ldots +e_r$, et par cons\'equent $g_i =
w_ig_1w_i^{-1}$; il en r\'esulte facilement que tous les $g_i$ appartiennent \`a la m\^eme
composante connexe du groupe $G(\Omega)$. Enfin, on v\'erifie ais\'ement que $\P(e_1 +
\ldots - e_i+ \ldots -e_j+ \ldots +e_r) = g_ig_j$, par cons\'equent les transformations $\P(e_1 +
\ldots - e_i+ \ldots -e_j+ \ldots +e_r)$ appartiennent \`a $G$.

\bigskip

\bigskip

\noindent $\underline {\hbox {B. La d\'ecomposition de Gauss.}}$

\bigskip

L'alg\`ebre de Lie $\goth G = \{ V \square V \}$ est r\'eductive d'involution de Cartan $\theta :X
\mapsto -X^*$, de d\'ecomposition de Cartan $\goth G = \goth k \oplus \goth p$ avec: 
$\goth k = \{ [\L(x),\L(y)] / x,y \in V \}$ et $\goth p = \{\L(x) / x\in
V\}$. Consid\'erons la d\'ecomposition de Cartan de $\goth G^1$: $\goth G^1 = \goth k \oplus \goth
p^1$ o\`u $\goth p^1 = \goth p \cap \goth G^1$. 

\noindent On pose:
$${\goth a}^1 = \{\L(\sum_{i=1}^r a_ie_i) /a_i \in \Bbb R , \sum_{i=1}^r a_i = 0 \}$$
c'est un sous-espace de Cartan de $\goth p^1$; les racines restreintes associ\'ees \'etant les
${\gamma}_{ij} = { {\epsilon}_i -{\epsilon}_j \over 2}, 1\leq i,j \leq r$ (avec
${\epsilon}_i(\L(e_i)) = \delta_{ij}$). On a une d\'ecomposition d'Iwasawa (associ\'ee
\`a un choix d'ordre sur les racines):  $\goth G^1 = \goth k \oplus \goth a^1 \oplus
\goth n$ avec:
$ \goth n = \sum_{i<j}\goth g_{ij}$ o\`u $\goth g_{ij} = V_{ij} \square e_j$. Rappelons que ce choix
n'est pas celui de [Fa.-Ko.]: notre $ \goth n$ correspond \`a leur $\bar {\goth n}$.

Soit $\goth h^-$ une sous-alg\`ebre ab\'elienne
maximale de $\goth m$ (o\`u $\goth m = \{X \in \goth k / \forall i, Xe_i = 0\}$), alors $\goth h
=
\goth h^-
\oplus
\goth a^1$ est une sous-alg\`ebre de Cartan de $\goth G^1$ et ${\goth h}_{\Bbb C}$ est une
sous-alg\`ebre de Cartan de ${\goth G^1}_{\Bbb C}$. On a une d\'ecomposition de Cartan de ${\goth
G^1}_{\Bbb C}$:
$${\goth G^1}_{\Bbb C} = (\goth k \oplus i\goth p^1)\oplus (i\goth k \oplus \goth p^1)$$ 
$i\goth h^- \oplus \goth a^1$ est un sous-espace de Cartan de $(i\goth k \oplus \goth p^1)$ et on a
la d\'ecomposition d'Iwasawa : ${\goth G^1}_{\Bbb C} = (\goth k \oplus i\goth p^1)\oplus
i\goth h^- \oplus \goth a^1 \oplus \goth n'$. Pour un certain choix de l'ordre sur les racines, on
aura: $\goth n' = {\goth n}_{\Bbb C} \oplus {\goth n}_{\goth m}$ o\`u $ {\goth n}_{\goth m}$ est
une sous-alg\`ebre de ${\goth m}_{\Bbb C}$ et o\`u ${\goth n}_{\Bbb C} = \sum_{i<j}V^{\Bbb C}_{ij}
\square e_j$.

Soient $N'$, $\bar N'$, $N_{\Bbb C}$, $N_{\goth m}$, $M_{\Bbb C}$ et $D$ les sous-groupes connexes
de
$G^1_{\Bbb C}$ d'alg\`ebres de Lie $\goth n'$, $\bar {\goth n'}$, ${\goth n}_{\Bbb C}$, ${\goth
n}_{\goth m}$, ${\goth m}_{\Bbb C}$ et ${\goth h}_{\Bbb C}$. Le groupe $M_{\Bbb C}$ normalise
$N_{\Bbb C}$ et  on a $N' =  N_{\Bbb C} .N_{\goth m} = N_{\goth m}.N_{\Bbb C}$.

\noindent La d\'ecomposition de Gauss de
$G^1_{\Bbb C}$ (voir [Ze.]) nous dit qu'il existe un ouvert $(G^1_{\Bbb C})_{\hbox {reg}}$ dense tel
que   
$$(G^1_{\Bbb C})_{\hbox {reg}} = {\bar N'}.D.N'.$$ 

\noindent Rappelons que tout \'el\'ement $ n$ de $N_{\Bbb C}$ s'\'ecrit de mani\`ere unique sous la
forme $n = \tau (z^{(2)}) \ldots  \tau (z^{(r)})$, avec $z^{(j)} \in
\bigoplus_{k=1}^{j-1}V^{\Bbb C}_{kj}$ et $\tau (z^{(j)}) = \exp(2z^{(j)} \square
e_j)$.

Enfin on sait ([Bl.2]) que les \'el\'ements $z$ de $V^{\Bbb C}$, inversible par rapport \`a tous
les
$\sum_{k=1}^ie_k$ ($1\leq i \leq r$), s'\'ecrivent sous la forme:
$$z = \tau (z_{(1)}) \ldots  \tau (z_{(r-1)})\sum_{k=1}^ra_ke_k$$
avec $\tau (z_{(j)}) = \exp(2z_{(j)}
\square e_j)$, o\`u $z_{(j)} \in \bigoplus_{k=j+1}^{r}V^{\Bbb C}_{jk}$, les $a_k$ \'etant des
nombres complexes non nuls; si l'\'el\'ement $z$ est dans $V$, les nombres
$a_k$ seront r\'eels et $z_{(j)} \in \bigoplus_{k=j+1}^{r}V_{jk}$.

Rappelons pour terminer, (voir par exemple [Ze.]) que si $\pi$ est une 
repr\'esenta-
\noindent tion holomorphe de poids dominant $\lambda$ de $G^1_{\Bbb C}$, on peut la voir comme
agissant dans un espace
$F_{\lambda}$ de fonctions analytiques sur
$N'$ de la mani\`ere suivante:
$\pi (g)f(z) = \lambda (zg)f(z_g)$ o\`u $f$ est un \'el\'ement de $F_{\lambda}$, $z$ est dans $N'$,
$g$ dans $(G)_{\Bbb C}$ (pourvu que $zg$ soit dans $(G^1_{\Bbb C})_{\hbox {reg}}$, et alors $z_g$
est la $N'$-composante de $zg$). Comme le
poids dominant est analytique, il est enti\`erement d\'efini par les valeurs qu'il prend sur $A^1$
et sur
$\exp (\goth h^-)$.

\bigskip

\bigskip

\eject

\noindent $\underline {\hbox {C. Les repr\'esentations du groupe $G$.}}$

\bigskip

Donnons pour commencer, la d\'ecomposition de l'espace ${\cal P}(V)$ des polyn\^omes sur
$V$, \`a valeurs complexes, sous l'action naturelle du groupe $G$. Les sous-espaces $V^{(k)} = \{
x
\in V / (e_1 + \ldots + e_k )x = x$ (resp. $V_{(k)}
 = \{ x \in V / (e_{r-k+1} + \ldots + e_r )x = x$) ($1 \leq k \leq
r$) sont des sous-alg\`ebres de $V$. Notons par $\Delta_k$ (resp. $\Delta^*_k$) la fonction polyn\^ome
sur $V$ d\'efinie par
$\Delta_k(x) = \det^{(k)}(x^{(k)})$ o\`u $x^{(k)}$ est la projection orthogonale de $x$ sur
$V^{(k)}$ et $\det^{(k)}$ est la fonction d\'eterminant de l'alg\`ebre $V^{(k)}$
(resp. $\Delta^*_k(x) = \det_{(k)}(x_{(k)})$ o\`u $x_{(k)}$ est la projection orthogonale de $x$
sur
$V_{(k)}$ et $\det_{(k)}$ est la fonction d\'eterminant de l'alg\`ebre $V_{(k)}$).

\noindent L'espace ${\cal P}(V)$ admet alors la d\'ecomposition en $G$-modules irr\'eductibles
(cf. [Fa.-Ko.]):

$${\cal P}(V) = \bigoplus_{\bf{m}=(m_1, \ldots ,m_r)}{\cal P}_{\bf m}$$
avec
${\cal P}_{\bf m} = \{\Delta _{\bf m}(gx) / g \in G \}$ (o\`u $\Delta _{\bf m}(x) =
{\Delta_1(x)}^{m_1
- m_2}{\Delta_2(x)}^{m_2 - m_3}...{\Delta_r(x)}^{m_r}$, ${\bf m} = (m_1,..,m_r)$ avec $ m_1\geq
m_2\geq ...\geq m_r \geq 0$ les $m_i$ \'etant des entiers). 

\noindent On notera $\pi _ {\bf m}$ l'action naturelle de $G$ dans l'espace ${\cal P}_{\bf m}$. Un
vecteur de plus haut poids pour le $G$ module ${\cal P}_{\bf m}$ est le polyn\^ome
$$\Delta^* _{\bf m} = {\Delta^*_1(x)}^{m_1
- m_2}{\Delta^*_2(x)}^{m_2 - m_3}...{\Delta^*_r(x)}^{m_r}$$
de poids dominant:
$${\lambda}_{\bf m}(\P (\exp a)) = a_1^{-2m_r}....a_r^{-2m_1}$$
($a = \sum a_ie_i$),
et le plus bas poids ${\mu}_{\bf m}$ par ${\mu}_{\bf m}(\P (\exp a)) = a_1^{-2m_1}....a_r^{-2m_r}$,
$\Delta _{\bf m}$ \'etant un vecteur de plus bas poids.

\noindent Remarquons que, consid\'er\'e comme
$G^1$ module, ${\cal P}_{\bf m}$ reste irr\'eductible et s'identifie \`a la repr\'esentation
${\cal P}_{\bf m^1}$ avec ${\bf m^1} = (m_1 - m_r,m_2 -
m_r,..,m_{r-1} - m_r,0)$ et que la d\'ecomposition $G =
{\Bbb R}^{+*}\times G^1$ nous donne ${\pi}_{\bf m} =  \chi ^{-{|{\bf m}|\over n}} \otimes
\pi^1 _ {\bf m^1}$, avec $\chi (g) = \Det g$. D'autre part, nous avons la

\proclaim Proposition II.1. La contragr\'ediente du $G^1$ module ${\cal P}_{{\bf m}}$
s'identifie au module ${\cal P}_{{\bf m^c}}$ o\`u ${\bf m^c} = (m_1 - m_r,m_1 - m_{r-1},..,m_1 -
m_2,0)$.

En effet, il n'est pas difficile de voir que la contragr\'ediente du $G$-module $({\cal P}_{{\bf
m}}, {\pi}_{\bf m})$ s'identifie \`a $({\cal P}_{{\bf m}}, \tau )$ avec $\tau (g) = {\pi}_{\bf
m}(g^{*-1})$. D'autre part, l'application $g \rightarrow g^*$ \'echange les groupes $N$ et $\bar
N$ et par cons\'equent un vecteur de plus haut poids pour $({\cal P}_{{\bf m}}, \tau )$  est
$\Delta _{\bf m}$; la proposition II.1 s'en d\'eduit ais\'ement. Notons qu'en utilisant la formule
(cf. [Fa.-Ko.], proposition VII.1.5) $\Delta _{\bf m}(x^{-1}) = \Delta^* _{-\bf m^ *}(x)$ (avec
${-\bf m^*} = (-m_r,
\ldots ,-m_1)$), on peut voir que l'application de 
${\cal P}_{{\bf m}}$ dans ${\cal P}_{{\bf m^c}}$:
$$P(x) \mapsto (\det x)^{m_1}P(x^{-1})$$
r\'ealise un isomorphisme de $G^1$ module entre la contragr\'ediente de ${\cal P}_{{\bf m}}$ et 
${\cal P}_{{\bf m^c}}$.

\bigskip

Soit maintenant $\pi$, une repr\'esentation irr\'eductible de dimension finie de $G$, elle
s'\'ecrit sous la forme $\chi ^{\alpha} \otimes \pi^1$ o\`u $\pi^1$ est une repr\'esentation
irr\'eductible du groupe $G^1$.

\proclaim Proposition II.2. La repr\'esentation  $\pi^1$ est la restriction d'une repr\'esentation
 holomorphe du groupe $G^1_{\Bbb C}$.

On utilise la classification. Le cas 5) est d\'emontr\'e dans [Kol.-Va.1]; le cas 6) r\'esulte du
fait qu'alors $G^1_{\Bbb C}$ est simplement connexe. Dans le cas 2), la repr\'esentation  $\pi^1$ de $\PSL
(r,{\Bbb R}) = SL (r,{\Bbb R})/ \{Id,(-1)^{r+1}Id \}$ peut \^etre consid\'er\'ee comme une repr\'esentation de $SL
(r,{\Bbb R})$ telle que $\pi^1 ((-1)^{r+1}g) = \pi^1 (g)$; sa d\'eriv\'ee
$d\pi ^1$ s'\'etend par lin\'earit\'e sur ${\goth G^1}_{\Bbb C}$, cette extension
s'int\`egre ( par simple connexit\'e) en une repr\'esentation holomorphe de $\SL (r,{\Bbb C})$  poss\'edant la
m\^eme parit\'e que $\pi^1$, elle redescend donc en une repr\'esentation de $G^1_{\Bbb C} = SL (r,{\Bbb C})/
\{Id,(-1)^{r+1}Id \}$; d\'emonstration
analogue pour les cas 3) et 4). 

\bigskip

Consid\'erons $\tilde \pi$ une repr\'esentation holomorphe de $G_{\Bbb C}$, elle est
caract\'eris\'ee par son poids dominant $\lambda _{\tilde \pi}$. Ecrivons $\lambda _{\tilde \pi}
(\P (\sum a_ie_i)) = \prod_1^r a_i^{m_i}$. Par analycit\'e de $\lambda _{\tilde \pi}$ les nombres
$m_i$ sont des entiers relatifs, et en \'ecrivant les conditions pour que $\lambda _{\tilde \pi}$
soit dominant, on trouve:
$m_1\leq m_2\leq \ldots \leq m_r$.

\proclaim Proposition II.3. Soit $\pi$ une repr\'esentation irr\'eductible de dimension finie de
$G$, sph\'erique (i.e. elle poss\`ede un vecteur non nul $K$ invariant). Alors elle peut s'\'ecrire
$$ \pi = \chi ^{\beta} \otimes \pi_{\bf m}.$$

Comme $\pi = \chi ^{\alpha} \otimes \pi^1$ est sph\'erique, la repr\'esentation $\pi^1$ l'est aussi; par la
proposition II.2 on peut se ramener au cas o\`u
$\pi^1$ est la restriction d'une repr\'esentation
${\tilde \pi}^1$ holomorphe de $G_{\Bbb C}$, ${\tilde \pi}^1$ est donc sph\'erique. En utilisant le th\'eor\`eme
de Cartan-Helgason on trouve que la restriction de 
$\lambda _{{\tilde \pi}^1}$ \`a $\goth h^-$ est nulle et que ${(\lambda _{{\tilde \pi}^1},\alpha) \over (\alpha,
\alpha)}
\in {\Bbb Z}_+$. Cette derni\`ere condition nous fournit ${m_j - m_i \over 2} \in {\Bbb Z}_+$ pour
$j>i$, c'est
\`a dire que les $m_j - m_i$ sont pairs, ce qui ach\`eve la d\'emonstration (voir aussi [Hi.-Ne.]
pour ce type de r\'esultat).

\bigskip

\bigskip

\noindent $\underline {\hbox {D. Les distributions z\^eta g\'en\'eralis\'ees.}}$

\bigskip

On note par $\Phi_j^{\bf m}(.,s)$ les distributions temp\'er\'ees
d\'efinies par prolongement m\'eromorphe des int\'egrales (qui sont convergentes pour $\Re s \geq
0$):

$$\int_{\Omega_j}{|\det x|^s \Delta _{\bf m}(x)f(x)dx} .$$
Pour ${\bf m} = {\bf 0}$, on retrouve les distributions z\^eta classiques, not\'ees
$\Phi_j(.,s)$ dans [Bl.1] et on a:
$$ \Phi_j^{\bf m}(.,s) = \Delta _{\bf m}\Phi_j(.,s)$$
Les p\^oles possibles de ces distributions sont ceux de la fonction 

$$\Gamma_{\Omega}(s+n/r+{\bf m}) = (2\pi)^{n-r\over2}\prod_{j=1}^{r} {\Gamma(s + m_j + 1 +
{(r-j)d\over 2})}$$ 
(voir [Fa.-Ko.] pour les propri\'et\'es de la fonction gamma du c\^one $\Omega$), et si $s_0$ est
p\^ole d'ordre
$o_{\bf m}(s_0)$ de cette fonction, il est facile de voir que sa multiplicit\'e pour  $\Phi_j^{\bf
m}(.,s)$, ne d\'epasse pas  $o_{\bf m}(s_0)$ cela r\'esulte de l'identit\'e de Bernstein (voir
[Fa.-Ko.], page 129):
$$(\det (\partial))^k|\det x|^{s+k}\Delta _{\bf m}(x) =(-1)^{j}|\det x|^s\Delta _{\bf m}(x){\Gamma
_{\Omega}(s+k+{\bf m}+{n\over r}) \over \Gamma
_{\Omega}(s+{\bf m}+{n\over r})}
\ ,\quad x\in
\Omega_j \  $$
o\`u, si $P$ est un polyn\^ome sur $V$, $P( \partial)$ est l'op\'erateur
diff\'erentiel \`a coefficients constants tel que: $P( \partial) \exp(<.,y>) = P(y)\exp(<.,y>)$.

Rapelons bri\`evement, pour la
commodit\'e du lecteur, la d\'efinition du changement de variable $\Phi$, d\'efini dans [Bl.1]:
on fixe un syst\`eme complet d'idempotents orthogonaux de $V$, $\{e_1,\ldots ,e_r\}$; dans la d\'ecomposition de
Peirce associ\'ee $V = \bigoplus _{1\leq i\leq j\leq r}V_{ij}$, l'alg\`ebre:
 $$\displaystyle V' = \{x\in V /e_1x = 0 \} = \bigoplus _{2\leq i\leq j\leq r}V_{ij}$$
est une sous-alg\`ebre de V simple, euclidienne et de rang $r-1$; posons enfin:
$\displaystyle W = \{x\in V /e_1x = {1\over 2}x\} = \bigoplus _{2\leq j\leq r}V_{1j}$;
l'application $\Phi$ est alors d\'efinie par:
$${\eqalign { \Phi : &{\Bbb
R}\times W\times V' \longrightarrow V\cr
 &(u,z,v) \longmapsto \exp(2z \square e_1)(ue_1
+v).\cr
}}$$
cette application r\'ealise
un $C^{\infty}$ diff\'eomorphisme de ${\Bbb R}^*\times W\times V'$ sur  
l'ouvert $ U = \{x \in V  / x_{11} \not = 0 \}$ (cf. la proposition 3.1 de [Bl.1]; pour un \'el\'ement $x$ de $V$, 
$x_{ij}$ d\'esigne la composante sur $V_{ij}$ de cet \'el\'ement suivant la d\'ecomposition de
Peirce).

\noindent Notons $\Phi_j^{\bf m}(.,s)_{|U}$ la restriction \`a $U$ de la distribution 
$\Phi_j^{\bf m}(.,s)$, et effectuons le changement de variable $\Phi$ dans l'ouvert $U$, on montre
alors (m\^eme d\'emonstration que celle du corollaire 3.2 de [Bl.1]) que l'image r\'eciproque par
$\Phi$ de $\Phi_j^{\bf m}(.,s)_{|U}$ not\'ee $\Phi^*[\Phi_j^{\bf m}(.,s)_{|U}] $ v\'erifie la
proposition:

\proclaim Proposition II.4. On a, pour $ (0<j<r)$:
$$\Phi^*[\Phi_j^{\bf m}(.,s)_{|U}] = x_{-}^{s+d(r-1)+m_1}\otimes
dz\otimes \Phi_{j-1}^{\bf m'}(.,s) +
x_{+}^{s+d(r-1)+m_1}\otimes dz\otimes
\Phi_j^{\bf m'}(.,s)$$
et
$$\Phi^*[\Phi_0^{\bf m}(.,s)_{|U}] = 
x_{+}^{s+d(r-1)+m_1}\otimes dz\otimes
\Phi_0^{\bf m'}(.,s)$$
$$\Phi^*[\Phi_r^{\bf m}(.,s)_{|U}] = x_{-}^{s+d(r-1)+m_1}\otimes
dz\otimes \Phi_r^{\bf m'}(.,s) $$
avec ${\bf m'} = (m_2, \ldots, m_r)$, et o\`u les distributions $\Phi_j^{\bf m'}$ sont 
d\'efinies sur l'alg\`ebre $V'$ comme les $\Phi_j^{\bf m}$ le sont sur $V$.

\noindent Dans cette proposition, nous avons utilis\'e les notations classiques $x_+^s$ (resp.
$x_-^s$) pour d\'esigner les distributions d\'efinies par

 $$x_+^s (\phi) = \int_0^{\infty}{x^s \phi(x) dx} \ \ \ \ ({\hbox {resp.}}\ \ x_-^s (\phi) =
\int_{-\infty}^0{\left|x\right|^s \phi(x) dx})$$

\noindent Cette proposition nous permettra, dans le paragraphe IV, de d\'eterminer, par une
r\'ecurrence sur le rang $r$, les p\^oles effectifs des distributions $\Phi_j^{\bf m}(.,s)$.

D\'efinissons la transform\'ee de
Fourier d'un \'el\'ement
$f$ de l'espace de Schwartz de $V$ par: ${\hat f} (y) =\int_V {e^{-i<x,y>}f(x)dx}$; \`a partir de 
l'\'equation fonctionnelle ([Sa.-Sh.], voir aussi [Sa.-Fa.]) des distributions $\Phi_j(.,s)$,
on obtient celle des $\Phi_j^{\bf m}(.,s)$ (je remercie J. Faraut de m'avoir signal\'e ce
r\'esultat):

\proclaim Proposition II.5. La transform\'ee de
Fourier de la distribution $\Phi_j^{\bf m}(.,s)$ est donn\'ee par :
$$\Phi_j^{\bf m}({\hat f},s-n/r) = \Gamma_{\Omega}(s+n/r+{\bf m})
\exp(-(rs+|{\bf m}|)i\pi/2)\sum_{k=0}^ru_{jk}(s)\Psi_k^{-\bf m^*}(f,-s)$$  o\`u les
fonctions
$u_{jk}(s)$ sont des polyn\^omes en $\exp(i{\pi}s)$, ${-\bf m^*} = (-m_r, \ldots ,-m_1)$ et avec 
$$\Psi_k^{{\bf m}}(f,s) = \int_{\Omega_j}{|\det x|^s \Delta^* _{\bf m}(x)f(x)dx}.$$

Dans  $\Phi_j({\hat f}\Delta _{\bf m},s-n/r)$, on commence par remplacer
${\hat f}\Delta _{\bf m}$ par $(-i)^{|{\bf m}|}(\Delta _{\bf m}(\partial)f{\hat)}$; 
l'\'equation fonctionnelle des distributions $\Phi_j(.,s)$ telle qu'\'enonc\'ee dans [Bl.1]
(th\'eor\`eme 2.1) conduit \`a la relation:

$$\Phi_j((\Delta _{\bf m}(\partial)f{\hat)},s-n/r) = \Gamma_{\Omega}(s) \exp(-{i\pi rs \over 2})
\sum_{k=0}^ru_{jk}(s)\Phi_k(\Delta _{\bf m}(\partial)f,-s).$$
Une int\'egration par parties et l'identit\'e de Bernstein:
$$\Delta _{\bf m}(\partial)|\det x|^s =  (-1)^{|{\bf m}|}{\Gamma_{\Omega}(-s+{\bf m}) \over
\Gamma_{\Omega}(-s)}\Delta _{\bf m}(x^{-1}) |\det x|^s$$
nous donne:
$$\Phi_k(\Delta _{\bf m}(\partial)f,-s) = {\Gamma_{\Omega}(s+{\bf m}) \over
\Gamma_{\Omega}(s)}\Phi_k(\Delta _{\bf m}(x^{-1})f(x), -s).$$
La proposition r\'esulte alors de la formule $\Delta _{\bf m}(x^{-1}) = \Delta^* _{-\bf m^*}(x)$
d\'eja \'evoqu\'ee plus haut.

Remarque II.6 On voit que dans cette \'equation fonctionnelle g\'en\'eralis\'ee,
les coefficients $u_{jk}(s)$ sont ceux de l'\'equation fonctionnelle des distributions
$\Phi_j(.,s)$, ce fait est essentiel pour les d\'emonstrations des r\'esultats du paragraphe IV.

\proclaim D\'efinition II.7. Soit $(\pi , L_{\pi})$ une repr\'esentation irr\'eductible de
dimension finie du groupe $G$. Une
 distribution vectorielle $T$ sur
$V$ \`a valeurs dans
$L_{\pi}$ est dite 
$G$-homog\`ene de degr\'e
$s$ si pour toute fonction test $\phi$,
$$<T,{\phi}_g> = (\Det g)^{{rs\over n}+1}\pi (g)<T,{\phi}>$$ o\`u
$g$ est dans $G$ et avec ${\phi}_g(x) = {\phi}(g^{-1}x)$. On notera ${\cal H}_{\pi}^s$ l'espace
de ces distributions. La distribution $T$ est dite $G^1$ invariante si 
$$<T,{\phi}_g> = \pi (g)<T,{\phi}>$$ o\`u
$g$ est dans $G^1$. On notera ${\cal I}_{\pi}$ l'espace
de ces distributions.

\noindent Dans le cas o\`u la repr\'esentation $\pi$ est triviale, on retrouve la notion habituelle
de distribution homog\`ene, on sait alors (voir par exemple [Bl.1]) d\'eterminer compl\`etement
l'espace ${\cal H}_{\pi}^s$ cela se fait gr\^ace \`a la consid\'eration des distributions
$\Phi_j(.,s)$. Dans ce travail, nous verrons
le r\^ole crucial jou\'e par les 
$\Phi_j^{\bf m}(.,s)$ pour la d\'etermination de ${\cal H}_{\pi}^s$ dans le cas g\'en\'eral.

Dans le paragraphe V. nous aurons besoin de la notion de distribution quasi-homog\`ene, introduite
dans [Ri.-St.]:

\proclaim D\'efinition II.8. Une distribution vectorielle $K$ sur $V$ \`a valeurs dans ${\cal
P}_{\bf m}$ est $G$ quasi-homog\`ene de degr\'e
$s$ s'il existe des distributions $U_1, \ldots ,U_l$ telles que pour tout $g$ dans $G$:
$$g.K = (\Det g)^{{rs\over n}+1}\pi_{\bf m}(g)K  + (\Det g)^{{rs\over n}+1} \sum_{i=1}^l(\Log \Det
g)^i\pi_{\bf m}(g)U_i.$$ 

\bigskip

\bigskip

\beginsection {III. Caract\'erisation des repr\'esentations.}

On se propose de caract\'eriser dans ce paragraphe les repr\'esentations $\pi$ pour lesquelles
l'espace ${\cal H}_{\pi}^s$ n'est pas trivial. On va retrouver des r\'esultats analogues \`a ceux
de [Kol.-Va.1].

\proclaim Lemme III.1. Si l'espace ${\cal H}^s_{\pi}$ contient une distribution non nulle de
support $\{0\}$, la
repr\'esentation $\pi$ est sph\'erique.

Soit $T$ une telle distribution. Dans une certaine base $e_{\alpha}$ de l'espace de la
repr\'esentation $\pi$, on a: $T = \sum P_{\alpha}(\partial) \delta e_{\alpha}$ o\`u les
$P_{\alpha}$ sont des polynomes. La repr\'esentation $\pi$ s'\'ecrit $\pi = \chi ^{\beta} \otimes
\pi ^1$ o\`u $\pi^1$ est une repr\'esentation  $G^1$. En \'ecrivant la condition
d'homog\'en\'eit\'e avec $g = \lambda Id$, on montre facilement que tous les polynomes $P_{\alpha}$
sont homog\`enes de m\^eme degr\'e ( on peut v\'erifier que ce degr\'e vaut $-n({rs\over n}+1 +
\beta )$). Ecrivons ensuite la condition
d'homog\'en\'eit\'e pour la distribution $T$ avec $g$ dans $G^1$, en notant $\lambda _{\alpha
\beta}(g)$ les coefficients matriciels de ${\pi}^1(g)$ dans la base $(e_{\alpha })$, il vient:
$$\sum (P_{\alpha}(\partial) \delta)(g.\phi) e_{\alpha} = \sum_{\alpha, \beta} \lambda _{\alpha
\beta}(g)(P_{\beta}(\partial) \delta )(\phi)
e_{\alpha}.$$
On en tire $$g.P_{\alpha}(\partial) = (\sum _{\beta}\lambda _{\alpha
\beta}(g^{-1})P_{\beta})(\partial)$$ o\`u $(g.P_{\alpha}(\partial))(\phi) =
g.(P_{\alpha}(\partial)(g^{-1}\phi))$. Mais $g.P_{\alpha}(\partial) =
(g^{*-1}P_{\alpha})(\partial)$ o\`u $g.P$ est l'action naturelle de $G$ sur les polynomes (voir
[Fa.-Ko.] page 222), tout ceci implique 
$$g^{*-1}P_{\alpha} = \sum _{\beta}\lambda _{\alpha \beta}(g^{-1})P_{\beta})$$
c'est \`a dire que sous l'action du groupe
$G^1$, les polynomes $P_{\alpha}$ engendrent la repr\'esentation  ${\pi}^1$;
il en r\'esulte  que cette derni\`ere est isomorphe \`a une repr\'esentation $\pi^1_{\bf m}$.

Dans la suite de l'article nous dirons que le support d'une
distribution est de rang $p$ s'il contient au moins un
\'el\'ement de rang $p$ mais aucun \'el\'ement de rang strictement
sup\'erieur \`a $p$.

\proclaim Lemme III.2. Supposons que l'espace ${\cal H}^s_{\pi}$ contienne une distribution non
nulle 
\`a  support de rang $p$ avec $1 \leq p \leq r-1$, alors la
repr\'esentation $\pi$ est sph\'erique.

La d\'emonstration consiste \`a adapter la preuve du lemme 5 de [Ri.-St.]. Soit $T$ une
distribution donn\'ee par le lemme. Dans un voisinage convenable $W$ du point $o_{p,q}$, la
distribution $T$ s'\'ecrit par le th\'eor\`eme de Schwartz comme somme finie:
$$T = \sum D_j{\Phi}_j$$
o\`u les ${\Phi}_j$ sont des distributions sur $S_{pq}$ et les $D_j$ sont des d\'eriv\'ees
partielles par rapport \`a des coordonn\'ees transverses \`a $S_{pq}$ en $o_{p,q}$.

On montre dans un premier temps que les distributions ${\Phi}_j$ sont \`a densit\'e $C^{\infty}$
sur $S_{pq} \cap W$. En effet, soit $\eta $ une fonction $C^{\infty}$ sur $G$, \`a support
compact convenable et telle que 
$$ B = \int_G (\Det g)^{{rs \over n}+1}\eta (g) \pi (g) dg $$
soit inversible. Utilisant l'homog\'en\'eit\'e de $T$, il vient facilement:
$$T(\phi) = B^{-1}\int_G \eta (g) T(g.\phi)dg$$
Pour toute fonction $\phi$ \`a support convenable. On conclut alors comme dans [Ri.-St.].

On arrive ainsi \`a la repr\'esentation suivante de $T$ dans $W$:

$$T = \sum_{m=0}^M {\sum_{i_1 \leq i_2 \leq \ldots \leq i_m}Y_{i_1}\ldots Y_{i_m} \psi _I}$$
 o\`u $I = (i_1,\ldots ,i_m)$,
les fonctions $\psi_I$ \'etant $C^{\infty}$ sur $S_{pq} \cap W$,
$Y_{i_1}\ldots Y_{i_m} $ sont des champs de vecteurs sur $W$ que l'on suppose engendrer en tout
point de $S_{pq} \cap W$ l'espace normal
\`a l'orbite $S_{pq}$, et o\`u $M$ est l'ordre de transversalit\'e de la distribution $T$.  

Exprimant l'homog\'en\'eit\'e de $T$ comme dans [Ri.-St.], nous pouvons \'ecrire pour
$\left|I\right| = M$ (cf. la relation (37) de [Ri.-St.]):
$$\pi (g) \psi_I(x) = (\Det g)^{-{rs \over n}-1} J_g(x) \sum_{\left|L\right|  = M}{\gamma _{LI}(g,x)
\psi_L(g.x)}  \leqno (1)$$
$g$ \'etant proche de
l'identit\'e, $x$ \'etant dans $S_{pq} \cap W$, o\`u $J_g(x)$ est le jacobien de l'action de $g$
sur
$S_{pq}$ et avec
$\gamma _{LI}(g,x) =
\beta _{IL}(g^{-1},x)$ les coefficients 
$\beta _{IL}(g,x)$ \'etant donn\'es par l'action de $g$ sur les champs de vecteurs:

$$g{\tilde Y}^I = \sum_{\left| J \right| = M}\beta _{JI}(g,x){\tilde Y}^J 
+ \sum_{\left| J \right| <M}{\ldots}$$
(${\tilde Y}^I$ d\'esigne l'adjoint de $Y^I$ sur $V$).

Dans la relation (1) faisons $x = o_{p,q}$, et prenons $g$ dans $H_{pq} = \{g \in G / go_{p,q} =
o_{p,q}\}$ et proche de l'identit\'e, on obtient:
$$\pi (g) \psi_I(o_{p,q}) = (\Det g)^{-{rs \over n}-1} J_g(o_{p,q}) \sum_{\left|L\right|  = M}{\gamma
_{LI}(g,o_{p,q}) \psi_L(o_{p,q})}.   \leqno (2)$$
Notons que, par analycit\'e, la relation  (2) est valable pour tout $g$ dans la composante neutre
$(H_{pq})_o$ du groupe $H_{pq}$; par ailleurs elle est \'egalement valable pour tout $g$ de
$H_{pq}$ qui envoie $S_{pq} \cap W$ dans lui-m\^eme.

Rappelons que nous avons not\'e $V_{(r-p)}$
l'espace normal \`a $S_{pq}$ en $o_{p,q}$, c'est  
une alg\`ebre de Jordan simple et euclidienne de rang
$r-p$. Dans la suite, si $O$ est un objet relatif \`a l'alg\`ebre $V$,
nous noterons par
$O_{(r-p)}$ l'objet correspondant pour l'alg\`ebre $V_{(r-p)}$; en particulier,
${\cal P}_{(r-p)}^{k}$ d\'esigne l'espace des polyn\^omes homog\`enes de degr\'e $k$ sur
$V_{(r-p)}$.

Dans une base convenable $\{ e_I,\left| I \right| = M\}$ de l'espace
${\cal P}_{(r-p)}^{M}$, l'application lin\'eaire de matrice $\gamma _{IJ}(g,o_{p,q})$, que nous
notons par 
$\tilde {\sigma}(g)$, n'est autre que la
projection sur l'espace ${\cal P}_{(r-p)}^{M}$ de l'action naturelle de $g$ sur ${\cal
P}^{M}$. En particulier, si l'on note par $H^{'}_{pq}$ le sous-groupe de $H_{pq}$ form\'e des
\'el\'ements de $H_{pq}$ laissant stable l'espace $N_{pq}$, et si l'on prend $g$ dans
$H^{'}_{pq}$, on a 
$\tilde {\sigma}(g) = \sigma (g_{(r-p)})$ o\`u $g_{(r-p)}$ est la restriction de $g$ \`a l'espace
$N_{pq} = V_{(r-p)}$ et o\`u $\sigma (g_{(r-p)})$ d\'esigne l'action naturelle de $g_{(r-p)}$ sur
${\cal P}_{(r-p)}^{M}$.

Consid\'erons l'application
lin\'eaire $L$ de ${\cal P}_{(r-p)}^{M}$ dans l'espace de la repr\'esenta-

\noindent tion $\pi$ telle que
$L(e_I) =
\psi_I(o_{p,q})$ par la relation (2) il vient imm\'ediatement la relation (3) pour $g$ dans
$(H^{'}_{pq})_o$: 
$$L \sigma (g_{(r-p)}) = (\Det g)^{{rs \over n}+1} J_g(o_{p,q})^{-1} \pi (g) L .  \leqno (3)$$
L\`a encore, la relation (3) est encore valable pour tout $g$ de
$H^{'}_{pq}$ qui envoie $S_{pq} \cap W$ dans lui-m\^eme.
Sous l'action de  $G_{(r-p)}$, ${\cal P}_{(r-p)}^{M}$ se d\'ecompose en modules irr\'eductibles:
$$ {\cal P}_{(r-p)}^{M} = \bigoplus _{\bf{n}=(n_{p+1}, \ldots ,n_r) \atop \left| \bf{n} \right| =
M} ({\cal P}_{(r-p)})_{\bf{n}}.$$

Observons que la restriction de $L$ \`a un espace $({\cal P}_{(r-p)})_{\bf{n}}$ est soit nulle,
soit injective; en effet supposons que
l'on ait $L(v) = 0$ pour un \'el\'ement $v$ non nul de cet espace, alors par la relation (3) il
vient
$L({\rm P}_{(r-p)}(x)v) = 0$ pour tous les
$x$ dans
$\Omega _{(r-p)}$ et donc, puisque l'application $x \rightarrow L({\rm P}_{(r-p)}(x)v)$ est
polynomiale, pour tous les
$x$ dans
$V^{\Bbb C}_{(r-p)}$; mais les \'el\'ements ${\rm P}_{(r-p)}(x)$ pour $x$ dans $V^{\Bbb
C}_{(r-p)}$ engendrent le groupe de structure de l'alg\`ebre
$V^{\Bbb C}_{(r-p)}$ et par cons\'equent $L$ est nulle sur l'espace $({\cal
P}_{(r-p)})_{\bf{n}}$ tout entier.

Supposons donc  que $L$ soit injective sur
l'espace $({\cal P}_{(r-p)})_{\bf{n}}$; la consid\'eration du vecteur 
$L(v_{\bf{n}}) = f$ o\`u $v_{\bf{n}}$ est un vecteur de plus haut poids de
$({\cal P}_{(r-p)})_{\bf{n}}$ va nous permettre de montrer que le poids dominant $\lambda$
de la repr\'esentation $\pi$ est trivial sur
$\exp (\goth h ^-)$.En effet suivant les rappels du paragraphe I, on
peut voir la repr\'esentation $\pi$ comme agissant dans un espace
$F_{\lambda}$ de fonctions analytiques sur $ N'$.
Un \'el\'ement 
$m$ de
$\exp (\goth h ^-)$ appartient \`a ${(H^{'}_{pq})}_o$ et le vecteur de plus haut poids $v_{\bf{n}}$
est invariant sous
$\sigma(m_{r-p})$ puisque $m_{(r-p)}$
appartient \`a
$M_{(r-p)}$. Par la relation (3) on aura par cons\'equent:
 
$$\lambda (zm)f(z_m) = f(z).$$ 
Prenons $ z = e$, et remarquons que la $ N'$-composante de $m$ est $e$; si $f(e)$ est
diff\'erent de $0$, on trouve la relation
$$\lambda (m) = 1$$
ce qui est bien l'\'egalit\'e souhait\'ee.

Nous sommes ramen\'e \`a montrer que $f(e)$ est non nul. Proc\'edant comme plus haut, on montre
que $f$ est $N_{\goth m}$ invariant et donc d\'etermin\'ee par ses valeurs
sur $ N_{\Bbb C}$ et, par analycit\'e, par ses valeurs sur $ N$.

D'apr\`es les rappels sur la d\'ecomposition de Gauss, $f$ peut 
s'\'ecrire, pour $z$ dans $N$: 
$$f(z) = f(\tau (z^{(2)}) \ldots \tau (z^{(r)})).$$
Pour $j>p$, l'\'el\'ement $\tau (z^{(j)})$ est dans $(H_{pq})_o$ (c'est une simple v\'erification
 utilisant le lemme VI.3.1 de [Fa.-Ko.]), et
$\tilde {\sigma} (\tau (z^{(j)}))v_{\bf{n}} = v_{\bf{n}}$ puisque l'on peut voir $v_{\bf{n}}$
comme un vecteur de plus haut poids pour la repr\'esentation ${\pi}_{(n_{p+1}, \ldots ,n_r,0,
\ldots ,0)}$; il en r\'esulte la relation: 
$$f(z) = f(\tau (z^{(2)}) \ldots \tau(z^{(p)})).$$
Pla\c cons nous alors dans l'alg\`ebre $V^{(p)} = \bigoplus _{0 < i\leq j\leq p}V_{ij}$ ( si $O$ 
est un objet relatif \`a l'alg\`ebre $V$, nous notons 
$O^{(p)}$ l'objet correspondant pour l'alg\`ebre $V^{(p)}$), et 
consid\'erons dans cette alg\`ebre l'\'el\'ement $\tau (z^{(2)})
\ldots \tau(z^{(p)})o_{p,q}$. G\'en\'eriquement (i.e en dehors de $\prod_1^p \Delta_k = 0$)
cet
\'el\'ement est inversible par rapport aux
$\sum_{k=1}^ie_k$ ($1 \leq i \leq p$) et sera donc, par les rappels du paragraphe II, de la forme
${ \bar n}ao_{p,q}$ avec ${\bar n}$ dans $ \bar {N}^{(p)}$ et $a$ dans $A^{(p)}$,
on peut voir ces deux \'el\'ements 
 dans $\bar N$ et $A$; finalement on peut \'ecrire $\tau (z^{(2)})
\ldots \tau(z^{(p)}) = {\bar n}ak$ avec $k$ agissant trivialement
sur $V_{(r-p)}$ et donc $\sigma (k_{(r-p)})v_{\bf{n}} = v_{\bf{n}}$. La relation (3) entraine
alors que dans un ouvert de $N^{(p)}$:
$$f(z) = \lambda (zk^{-1})f(e)$$
car $zk^{-1} \in \bar {N'}.D$ et donc sa $ N'$-composante est
r\'eduite \`a e.
Il r\'esulte de tout ceci que si $f$ est non nulle sa valeur en $e$ l'est aussi. Noter que cette
d\'emonstration donne l'unicit\'e de $f$ \`a un scalaire multiplicatif pr\`es

Ecrivons maintenant $\lambda (\P (\sum a_ke_k)) = \prod_{k=1}^r{a_k^{m_k}}$, o\`u les $m_k$ sont
des entiers relatifs tels que $m_1 \leq m_2 \leq \ldots \leq m_r$, il s'agit de voir que les
entiers relatifs
$m_k$ sont de m\^eme parit\'e. Pour cela, soit $g_i = \P(-e_1 + \ldots -e_i + \ldots + e_r)$,
c'est un \'el\'ement du groupe $H^{'}_{pq}$ qui envoie $S_{pq} \cap W$ dans lui-m\^eme; d'autre
part $(g_i)_{(r-p)} = {\rm P}_{(r-p)}(e_{p+1} + \ldots -e_i + \ldots + e_r)$ (ou $Id$ si $p>i$).
Par la relation (3) il vient:
$$\sigma (\rm {P}_{(r-p)}(e_{p+1} + \ldots -e_i + \ldots + e_r))v_{\bf{n}} = v_{\bf{n}}$$
puisque $v_{\bf{n}}$ est un vecteur de plus haut poids de $({\cal P}_{(r-p)})_{\bf{n}}$. Par
cons\'equent:
$$f = \pi (g_i) f.$$
En regardant en e, on arrive \`a:
$$\lambda ( {\rm P}_{(r-p)}(e_{p+1} + \ldots -e_i + \ldots + e_r)) = 1.$$
Ce qui entraine que $m_1$ et $m_i$ ont m\^eme parit\'e et ach\`eve la
d\'emonstration du lemme.

\bigskip

\bigskip

Nous sommes maintenant en mesure de montrer le r\'esultat essentiel de ce paragraphe.

\proclaim Th\'eor\`eme III.3. L'espace ${\cal H}^s_{\pi}$ est non nul si et seulement si la 
 repr\'esentation $\pi$ est sph\'erique.

Soit $T$ un \'el\'ement non nul de
${\cal H}^s_{\pi}$, par les lemmes III.1 et III.2, nous pouvons supposer que la restriction $T_i$
de la distribution $T$  \`a l'une des orbites ouvertes $\Omega _{i}$ est non nulle. Par un
r\'esultat classique (voir par exemple [Ra.-Sc.], lemme 5.12), la repr\'esentation
${\chi}^{{rs\over n}+1 +
\alpha }
\bigotimes {\pi}^1$ poss\`ede alors un vecteur $K_i$-invariant non nul, o\`u l'on a pos\'e
$K_i = \{g \in G / go_{r-i,i} = o_{r-i,i} \}$. Le lemme 5.12 de [Ra.-Sc.] montre aussi que
le vecteur $K_i$ invariant d\'etermine d'une mani\`ere unique la distribution $T_i$. Comme
le groupe
$K_i$ est inclus dans
$G^1$, il en r\'esulte que la repr\'esentation ${\pi}^1$ poss\`ede un vecteur $a_i$ non nul et
$K_i$-invariant. On sait par la proposition II.2 que ${\pi}^1$ s'\'etend en une repr\'esentation
holomorphe du groupe $G^1_{\Bbb C}$  encore not\'ee ${\pi}^1$. Cette repr\'esentation poss\`ede
donc un vecteur $(K_i)_{\Bbb C}$-invariant, o\`u l'on a not\'e par $(K_i)_{\Bbb C}$ le sous groupe
de
$G^1_{\Bbb C}$ correspondant \`a l'alg\`ebre de Lie $({\goth k}_i)_{\Bbb C}$, complexifi\'ee de
l'alg\`ebre de Lie ${\goth k}_i$ de $K_i$. Il suffit alors de remarquer que $K_{\Bbb C}$ et 
$(K_i)_{\Bbb C}$ sont conjugu\'es dans $G^1_{\Bbb C}$ pour pouvoir conclure que la
repr\'esentation $\pi^1$ poss\`ede un vecteur $K$-invariant non nul. Notons qu'ici aussi nous avons
l'unicit\'e de
$T_i$ \`a un scalaire multiplicatif pr\`es, puisque il en est de m\^eme du vecteur $K$-invariant.

Remarque III.4: Il n'est pas difficile d'adapter les d\'emonstrations effectu\'ees pour
les distributions $G^1$ invariantes, nous obtenons ainsi la g\'en\'eralisation du r\'esultat de
[Kol.-Va.1]:

\proclaim Th\'eor\`eme III.5. L'espace ${\cal I}_{\pi}$ est non nul si et seulement si la 
 repr\'esentation $\pi$ est sph\'erique.

\bigskip
\bigskip

\beginsection {IV. Construction de distributions homog\`enes.}

Notons par ${\cal H}_{\bf m}^s$ l'espace des distributions $G$-homog\`enes de degr\'e $s$ \`a
valeurs dans ${\cal P}_{\bf m}$. Le but de ce paragraphe est de construire des \'el\'ements de
${\cal H}_{\bf m}^s$. La strat\'egie est analogue \`a celle utilis\'ee dans [Ri.-St]. On commence
par une \'etude pr\'eliminaire.

\noindent $\underline {\hbox {A. Une \'etude scalaire pr\'eliminaire.}}$

On consid\`ere une distribution de la forme:
$$A^{\bf m}(.,s) = \sum_{j=0}^{r}c_j\Phi_j^{\bf m}(.,s).$$
Pour un p\^ole effectif $s_0$ d'ordre $\alpha$ de
la  distribution $A^{\bf m}(.,s)$, notons  
$$A^{\bf m}(.,s) = \sum_{h=1}^{h=\alpha}{A^{\bf
m}_h(.,s_0)\over(s-s_0)^h} + A^{\bf m}_0(.,s_0) + \ldots $$ 
son d\'eveloppement de
Laurent au voisinage de ce p\^ole. Nous voulons connaitre les p\^oles effectifs de la
distribution $A^{\bf m}(.,s)$ et avoir des pr\'ecisions sur le support des
distributions $A^{\bf m}_h(.,s_0)$ (en consid\'erant la restriction de $A^{\bf m}(.,s)$ \`a un
ouvert $\Omega_j$, on voit que le support de la distribution $A^{\bf
m}_0(.,s_0)$ est de rang $r$).

Supposons d'abord l'entier $d$
pair et construisons le polyn\^ome $P$ satisfaisant aux deux
conditions suivantes (o\`u $d^{\circ}P$ d\'esigne le degr\'e de $P$):
$$\left\{ \matrix{P(j) = c_j\exp (i\pi s_0j) \cr 
\d^{\circ}P \leq  r \cr}\right.
$$

\noindent Rappelons d'autre part que pour un p\^ole $s_0$ de la fonction
$\Gamma_{\Omega}(s+n/r+{\bf m})$, nous avons not\'e $o_{\bf m}(s_0)$ sa multiplicit\'e.

\proclaim Th\'eor\`eme IV.1. La distribution $A^{\bf m}(.,s)$ poss\`ede en
$s_0$ un p\^ole d'ordre $$\min(\d^{\circ}P,o_{\bf m}(s_0)).$$

La d\'emonstration  \'etant identique \`a celle du th\'eor\`eme 3.3 de [Bl.1], nous ne donnons pas
les d\'etails pour lesquels nous renvoyons \`a l'article cit\'e. On se place d'abord dans le cas
o\`u l'on a $o_{\bf m}(s_0) = r$ c'est \`a dire que $s_0 \leq -m_1-n/r$, et l'on d\'emontre le
r\'esultat pour une base particuli\`ere $\{P_0, \ldots ,P_r \}$ des polyn\^omes de
degr\'e inf\'erieur ou \'egal \`a  $r$; dans cette \'etape, le seul changement par rapport \`a la
d\'emonstration du th\'eor\`eme 3.3 de [Bl.1] est l'utilisation de l'\'equation fonctionnelle
g\'en\'eralis\'ee (c'est-\`a dire la proposition II.5) pour les distributions $\Phi_j^{\bf m}(.,s)$
au lieu de l'\'equation fonctionnelle classique; c'est ici que la remarque II.6 est cruciale,
puisque la d\'emonstration du th\'eor\`eme 3.3 de [Bl.1] est essentiellement bas\'ee sur des
manipulations de certains des coefficients $u_{jk}(s)$ intervenant dans l'\'equation fonctionnelle
classique.

\noindent Si maintenant l'on a 
$$-n/r-m_1 < s_0 \leq -m_r$$
on consid\`ere comme dans [Bl.1], la restriction de la distribution $A^{\bf m}(.,s)$ \`a l'ouvert
$U = \{ x \in V / \  \ x_{11} \not = 0 \}$ et l'on montre (par une r\'ecurrence sur $r$), comme dans
le cas scalaire (en utilisant notre proposition II.4 au lieu du corollaire 3.2 de [Bl.1]) que cette
restriction pr\'esente en
$s_0$ un p\^ole d'ordre
$\min(\d^{\circ}P,o_{\bf m}(s_0))$; la distribution $A^{\bf m}(.,s)$ poss\`ede donc en $s_0$ un
p\^ole d'ordre au moins
\'egal \`a
$\min(\d^{\circ}P,o_{\bf m}(s_0))$ mais par ailleurs la multiplicit\'e ne peut pas d\'epasser ce
nombre puisqu'elle est major\'ee par $o_{\bf m}(s_0)$ et que  $A^{\bf m}(.,s) = \Delta _{\bf
m}(x)A^{\bf 0}(.,s)$, cela ach\`eve la d\'emonstration.

Supposons l'entier $d$ impair et construisons les deux polyn\^omes
$P_o$ et $P_1$ tels que :
 $$\left\{\matrix{c_j\exp (i{\pi}js_o) = P_o(j)
\qquad \hbox {si j est pair} \cr 
\d^{\circ}P_o \leq  [{r\over2}]\cr
c_j\exp (i{\pi}js_o) = P_1(j) \qquad \hbox{si j est impair} \cr
\d^{\circ}P_1 \leq  [{r-1\over 2}] \cr}\right.$$
Posons par ailleurs 
$\epsilon = \epsilon (P_0,P_1) = 1 \quad \hbox {si} \quad \d^{\circ} (P_o-P_1) =
\sup(\d^{\circ}P_o,\d^{\circ}P_1)$, et
$\epsilon = 0 $ sinon.

\proclaim Th\'eor\`eme IV.1'.
 1) Si $s_0$ est un demi-entier, la distibution $A^{\bf m}(.,s)$ y poss\`ede
un p\^ole d'ordre 
$$ \min(\sup(\d^{\circ} P_o,\d^{\circ} P_1) , o_{\bf m}(s_0))$$
 2) Si $s_0$ est un entier, la distibution $A^{\bf m}(.,s)$ y poss\`ede
un p\^ole d'ordre 
$$\min(\sup(\d^{\circ} P_o,\d^{\circ} P_1)+\epsilon , o_{\bf m}(s_0)).$$

 La d\'emonstration est identique au cas scalaire (th\'eor\`eme 5.1 de [Bl.1]), 
en employant la remarque II.6 et la proposition II.4.

En ce qui concerne le support des $A^{\bf m}_h(.,s_0)$, les r\'esultats suivants se d\'emontrent
\'egalement en copiant le cas scalaire:

\proclaim Th\'eor\`eme IV.2. Supposons l'entier $d$ pair. Pour tout p\^ole $s_0$ d'ordre $\alpha$
de
$A^{\bf m}(.,s)$, et tout entier
$h$, compris entre $1$ et
$\alpha$, le support de la distribution
$A^{\bf m}_h(.,s_0)$ est de rang $r-h$.
Supposons l'entier $d$ impair. Pour tout p\^ole $s_0$ d'ordre $\alpha$
de $A^{\bf m}(.,s)$, et tout entier $h$ compris entre $1$ et
$\alpha$, le support de la distribution $A^{\bf m}_h(.,s_0)$ est de rang
$r+1-2h$ si $s_0$ est un entier, et de rang $r-2h$ si $s_0$ est un demi-entier.

\bigskip

\noindent $\underline {\hbox {B. Construction de distributions homog\`enes \`a valeurs dans ${\cal
P}_{\bf m}$.}}$

Avec $(e_{\alpha})$ une base de ${\cal P}_{\bf m}$, et $(f_{\alpha})$ sa base duale dans ${\cal
P}_{{\bf m^c}}$ (par la Proposition II.1 on a identifi\'e le $G^1$ module ${\cal P}_{{\bf m^c}}$
\`a la contragr\'ediente de ${\cal P}_{\bf m}$) nous pouvons d\'efinir l'application
$h_{\bf m}$ de $V$ dans
${\cal P}_{\bf m}$ par:
$$h_{\bf m}(x) = \sum_{\alpha} f_{\alpha}(x)e_{\alpha}.$$
Par construction, cette application v\'erifie le lemme suivant:

\proclaim Lemme IV.3. L'application $h_{\bf m}$ est $G^1$ \'equivariante, d'une mani\`ere pr\'ecise:
$$h_{\bf m}(g.x) = \chi^{rm_1 \over n}(g)\pi_{\bf m}(g)h_{\bf m}(x)$$
pour $g$ dans $G$.

On consid\`ere alors les distributions sur $V$ \`a valeurs dans ${\cal P}_{\bf m}$ du
type:
$$T_j^{\bf m}(\phi,s) = \int_{\Omega_j}{|\det x|^{s-{m_1}} h_{\bf m}(x)\phi
(x)dx}.$$ 
Les p\^oles de ces distributions sont parmi 
 ceux des distributions 

\noindent $\Phi_j^{\bf m^c}(.,s-m_1)$, (${\bf m^c} = (m_1 - m_r,m_1 - m_{r-1},..,m_1 -
m_2,0)$) cela r\'esulte du fait que les polyn\^omes $(f_{\alpha})$ sont, pour un choix
convenable de la base $(e_{\alpha})$, de la forme $\Delta _{\bf m^c}(gx)$. Appelons un nombre
complexe
$s$ critique pour $\pi_{\bf m}$ s'il est un p\^ole de la fonction
$$\Gamma_{\Omega}(s+n/r-{\bf m^*}) = (2\pi )^{n-r \over 2}\prod^r_{j=1} \Gamma
(s+1-m_j+(j-1)d/2).$$
Par les th\'eor\`emes IV.1 et IV.1', un nombre complexe
$s$ est critique pour $\pi_{\bf m}$ s'il est un p\^ole possible pour les distributions $\Phi_j^{\bf
m^c}(.,s-m_1)$; on a donc:

\proclaim Proposition IV.4. La distribution $T_j^{\bf m}(.,s)$ se prolonge en tout point $s$
non critique en un \'el\'ement de l'espace ${\cal H}_{\bf m}^s$.

Si $s_0$ est p\^ole effectif $s_0$ d'ordre $\alpha$ de
 $T_j^{\bf m}(.,s)$, on v\'erifie que les distributions $T^{\bf
m}_{j \alpha}(.,s_0)$ sont $s_0$-homog\`enes, et que les autres coefficients $T^{\bf
m}_{jh}(.,s_0)$ sont quasi-homog\`enes de degr\'e $s_0$.

Donnons, pour terminer ce paragraphe, un r\'esultat utile pour la suite; supposons d'abord l'entier
$d$ pair.

\proclaim Proposition IV.5. Soient $1 \leq p \leq r-1$ et $0 \leq q \leq p$. Pour tout $s_0$
critique  pour $\pi_{\bf m}$ tel que  $ r-p \leq o_{\bf -m^*}(s_0)$, il existe une distribution
$G$-homog\`ene de degr\'e
$s_0$, de support de rang $p$ et contenant l'orbite $S_{p,q}$.

On consid\`ere une distribution du type
$$A^{\bf m}(.,s) = \sum_{j=0}^{r}(-1)^{sj}j^{r-p}T_j^{\bf m}(.,s)$$
telle que par un choix convenable de la base $(e_{\alpha})$, une de ses composantes soit
$\sum_{j=0}^{r}(-1)^{sj}j^{r-p}\Phi_j^{\bf m^c}(.,s-m_1)$, elle pr\'esentera en $s_0$ un p\^ole
d'ordre $r-p$ d'apr\`es le th\'eor\`eme IV.1, et par le th\'eor\`eme
IV.2, la distribution
$A_{r-p}^{\bf m}(.,s_0)$ est \`a support de rang
$p$. Pour la composante $\sum_{j=0}^{r}(-1)^{sj}j^{r-p}\Phi_j^{\bf m^c}(.,s-m_1)$ de la
distribution $A^{\bf m}(.,s)$, on peut copier la d\'emonstration du corollaire 3.7 de [Bl.1] ce qui
nous donne la proposition.

Dans le cas o\`u $d$ est impair, on a

\proclaim Proposition IV.5'. Soient $1 \leq p \leq r-1$ et $0 \leq q \leq p$. Supposons $p$ de la
forme $r+1-2h$ avec $h$ un entier; alors pour tout $s_0$ entier et critique 
pour $\pi_{\bf m}$ tel que  $ h \leq o_{\bf -m^*}(s_0)$, il existe une distribution $G$-homog\`ene
de degr\'e
$s_0$, de support de rang $p$ et contenant l'orbite $S_{p,q}$. Supposons $p$ de la forme $r-2h$ avec
$h$ un entier; alors, pour tout demi-entier $s_0$ tel que  $h \leq o_{\bf -m^*}(s_0)$, il existe une
distribution $G$-homog\`ene de degr\'e $s_0$, de support de rang $p$ et contenant l'orbite
$S_{p,q}$.

Supposons $p$ de la forme $r+1-2h$ avec $h$ un entier et soit $s_0$ un entier critique 
pour $\pi_{\bf m}$ tel que  $ h \leq o_{\bf -m^*}(s_0)$ ; 
on consid\`ere les distributions
$$A^{\bf m}(.,s) = \sum_{j \equiv 0 [2]}\exp (-i\pi s_0)j^{h-1}T_j^{\bf m}(.,s)$$ et 
$$B^{\bf m}(.,s) = \sum_{j \equiv 1 [2]}\exp (-i\pi s_0)j^{h-1}T_j^{\bf m}(.,s).$$ Par le
th\'eor\`eme IV.1', ces distributions pr\'esentent en $s_0 $ un p\^ole d'ordre $h$ puisque $ h \leq
o_{\bf -m^*}(s_0)$, et par le th\'eor\`eme IV.2, les distributions $A^{\bf m}_{h}(.,s_0)$ et
$B^{\bf m}_{h}(.,s_0)$ sont \`a support de rang
$r+1-2h = p$. On montre comme dans la proposition IV.5 que le support de  ces distributions contient
toutes les orbites de
rang $p$.

Dans le cas o\`u $p$ est de la forme $r-2h$ avec $h$ un entier, on
consid\`ere les distributions
$$A^{\bf m}(.,s) = \sum_{j \equiv 0 [2]}\exp (-i\pi s_0)j^{h}T_j^{\bf m}(.,s)$$ et 
$$B^{\bf m}(.,s) = \sum_{j \equiv 1 [2]}\exp (-i\pi s_0)j^{h}T_j^{\bf m}(.,s).$$
A nouveau, par le
th\'eor\`eme IV.1', ces distributions pr\'esentent en $s_0 $ un p\^ole d'ordre $h$ et par le
th\'eor\`eme IV.2, les distributions $A^{\bf m}_{h}(.,s_0)$ et $B^{\bf m}_{h}(.,s_0)$ sont \`a
support de rang
$r-2h = p$, les m\^emes arguments s'appliquent, mais cette fois \c ci on trouve que le support de
$A^{\bf m}_{h}(.,s_0)$ contient toutes les orbites
$S_{p,q}$ avec $q$ pair, et que celui de $B^{\bf m}_{h}(.,s_0)$ contient toutes les orbites
$S_{p,q}$ avec $q$ impair.

\bigskip

\bigskip

\beginsection {V. L'espace ${\cal H}_{\bf m}^s$ et les distributions invariantes \`a support
singulier.}

\bigskip

Le but de ce paragraphe est de montrer que la dimension de l'espace ${\cal
H}_{\bf m}^s$ est $r+1$.

\proclaim Lemme V.1.  Soit $T$ une distribution non nulle, $G$ quasi-homog\`ene de degr\'e $s$,
\`a valeurs dans ${\cal P}_{\bf m}$ et de support de rang $0<p<r$. Le nombre $s$ est alors
n\'ecessairement critique pour la repr\'esentation $\pi_{\bf m}$, avec 
$$s \leq -(n/r-dp/2) + m_{r-p}.$$
De plus, si $S_{pq}$ est une orbite de rang $p$, la restriction de $T$ 
\`a $W_{pq}$ est homog\`ene et unique \`a une constante multiplicative pr\`es.

La d\'emonstration consiste \`a reprendre celle du lemme III.2 et \`a suivre la d\'emarche de
Ricci et Stein. Supposons d'abord que
$T$ soit homog\`ene.  On v\'erifie facilement que les \'el\'ements de $A$ de la forme $ a =
\P(o_{p,0} +
\sum_{i = p+1}^ra_ie_i)$ ($a_i$ positifs non nuls) appartiennent \`a
$(H'_{pq})_o$, et que donc la relation (3) du lemme III.2 est valable pour ces \'el\'ements.
D'autre part, il est facile d'obtenir les formules 
$$\eqalign{
{\P (a)}_{(r-p)} = & {\P}_{(r-p)}(\sum_{i = p+1}^ra_ie_i) \cr
\Det {\P}(a) = & (a_{p+1}\ldots a_r)^{2n \over r} \cr
J_{\P (a)}(o_{p,q}) = & (a_{p+1}\ldots a_r)^{dp}. \cr
}$$
Maintenant on a: ${\P (a)}_{(r-p)} v_{\bf n} = a_{p+1}^{-2n_r} a_{p+2}^{-2n_{r-1}} \ldots 
a_r^{-2n_{p+1}}v_{\bf n}$, puisque $ v_{\bf n}$ est un vecteur de plus haut poids de l'espace
$({\cal P}_{(r-p)})_{\bf{n}}$. La relation (3) du lemme III.2 nous donne alors en $z = e$ les
relations:
$$\matrix{
-2n_r & = & 2s-dp+2n/r-2m_{r-p} \cr
-2n_{r-1} & = & 2s-dp+2n/r-2m_{r-p-1} \cr
\ \vdots  &   &                 \cr
-2n_{p+1} & = & 2s-dp+2n/r-2m_1 \cr
} \leqno (I) $$
qui montrent que $s \leq -(n/r-dp/2) + m_{r-p}$, et que la donn\'ee de $s$, $p$, et
$\pi_{\bf m}$ d\'etermine d'une mani\`ere unique $\bf n$; la d\'emonstration du lemme III.2
prouve alors que la famille $(\psi_I(o_{p,q}))_{|I|=M}$ est unique \`a un facteur multiplicatif 
pr\`es; par la relation (1) de la d\'emonstration du lemme III.2 il en de m\^eme de la famille
$(\psi_I)_{|I|=M}$ sur
$S_{pq}\cap W$. Maintenant, par les propositions IV.5 et IV.5', on sait qu'il existe une
distribution
$U$, homog\`ene de degr\'e $s$ et de support de rang $p$ contenant $S_{pq}$, le r\'esultat
pr\'ec\'edent nous montre qu'il existe un scalaire $a$ telle que la distribution $T-aU$ ne
pr\'esente plus dans
$S_{pq}\cap W$ de termes d'ordre $M$, ce qui entraine que $T = aU$ dans $S_{pq}\cap W$ et par
homog\'en\'eit\'e dans $W_{pq}$ tout entier.

Pla\c cons nous maintenant dans le cas g\'en\'eral o\`u $T$ v\'erifie la relation:
$$g.T = (\Det g)^{{rs\over n}+1}\pi_{\bf m}(g)T  + (\Det g)^{{rs\over n}+1} \sum_{i=1}^l(\Log \Det
g)^i\pi_{\bf m}(g)U_i$$
et suivons encore Ricci et Stein en reprenant la d\'emonstration du lemme III.2. On montre
facilement que dans $W$:
$$T = \sum_{m=0}^M {\sum_{i_1 \leq i_2 \leq \ldots \leq i_m}Y_{i_1}\ldots Y_{i_m} \psi _I}$$
ainsi que: 
$$U_j = \sum_{m=0}^M {\sum_{i_1 \leq i_2 \leq \ldots \leq i_m}Y_{i_1}\ldots Y_{i_m} \mu _{Ij}}$$ 
o\`u $\psi _I$ et $\mu _{Ij}$ sont des fonctions r\'eguli\`eres sur $S_{pq}\cap W$. Comme dans
[Ri.-St], on montre ensuite que $\mu _{Ij}(o_{pq}) = 0$ ce qui nous permet de proc\'eder comme
dans le cas homog\`ene.

\bigskip

\proclaim Lemme V.2. Soit $T$ une distribution non nulle, $G$ homog\`ene de degr\'e $s$,
\`a valeurs dans ${\cal P}_{\bf m}$, et port\'ee par l'origine. Le nombre $s$ est alors
n\'ecessairement critique pour la repr\'esentation $\pi_{\bf m}$, avec 
$$s \leq -n/r + m_r.$$
De plus, une telle distribution est unique \`a une constante multiplicative pr\`es.

Soit en effet une distribution $T = \sum P_{\alpha}(\partial) \delta
e_{\alpha}$ homog\`ene de degr\'e $s$ et \`a valeurs dans ${\cal P}_{\bf m}$; la
d\'emonstration du lemme III.1 nous dit que tous les polyn\^omes $P_{\alpha}$ sont homog\`enes de degr\'e
$$l = -rs - n - |\bf {m}|.$$
et que le $G^1$ module engendr\'e par eux est $\pi_{\bf m^1}$. Le $G$ module engendr\'e par les $P_{\alpha}$ est
donc de la forme $(\det x)^a{\cal P}_{\bf m^1}$, o\`u $a$ est un entier positif tel que $l = ra +
|{\bf m}|-rm_r$; il en r\'esulte que:
$$s = -n/r + m_r -a$$
et par cons\'equent $s$ est critique pour la repr\'esentation $\pi_{\bf m}$. Si maintenant $(e_{\alpha})$ est une
base orthonorm\'ee pour le produit de Fischer (voir par exemple [Fa.-Ko.] pour le produit de Fischer), la
d\'emonstration du lemme III.1 nous montre \'egalement que: 
$$(\det x)^{-m_r} e_{\alpha} \mapsto (\det x)^{-a} P_{\alpha}$$
est $G^1$ \'equivariante, le lemme de Schur nous donne donc le
r\'esultat d'unicit\'e.

\proclaim Th\'eor\`eme V.3. La dimension de l'espace ${\cal H}_{\bf m}^s$ est
$r+1$.

La d\'emonstration est identique \`a celle du cas scalaire donn\'ee dans
[Ri-St.]. Indiquons en bri\`evement les grandes lignes dans le cas o\`u
l'entier $d$ est pair.

Soit donc $T$ un \'el\'ement non nul de l'espace ${\cal H}_{\bf m}^s$ et
supposons d'abord que $s$ ne soit pas critique pour la repr\'esentation $\pi
_{\bf m}$. Il r\'esulte du lemme V.1 que la restriction de $T$ \`a $V^{\times}$
est non nulle; par la d\'emonstration du Th\'eor\`eme III.3, il existe des
nombres $a_1, \ldots, a_r$ tels que sur $V^{\times}$:
$$T = \sum a_iT_i^{\bf m}(,s).$$
$s$ n'\'etant pas critique, on peut consid\'erer la distribution $T-\sum
a_iT_i^{\bf m}(,s)$ sur $V$ tout entier; elle est port\'ee par $S$ et
homog\`ene et donc par le lemme V.1, $s$ devrait etre critique; par cons\'equent on a sur $V$:
$$T = \sum a_iT_i^{\bf m}(,s).$$

On suppose maintenant que $s=s_0$ est critique et que la restriction de $T$ \`a $V^{\times}$
soit non nulle. Consid\'erons les 
distributions ($i=0,1,
\ldots ,r$):
$$A_i^{\bf m}(.,s) = \sum_{j=0}^{r}(-1)^{s_0j}j^iT_j^{\bf m}(.,s).$$
Le th\'eor\`eme IV.1 nous dit que la distribution $A_i^{\bf m}(.,s) $ pr\'esente en $s_0$
un p\^ole d'ordre $\min (i, o_{\bf m}(s_0))$. Ecrivons les d\'eveloppements de Laurent:

$$\matrix{
 A_0^{\bf m} & = & {} & {} &{}                   &+ A^{\bf m}_{00} & + \ldots \cr
 A_1^{\bf m} & = & {} & {} &{A^{\bf m}_{11} \over s-s_0} &+ A^{\bf m}_{10} & + \ldots \cr
 \ \vdots & & &    &                     &         &           \cr
 A_h^{\bf m} & = &{A^{\bf m}_{hh} \over (s-s_0)^h} & +\ldots + & {A^{\bf m}_{h1} \over s-s_0} & +
A^{\bf m}_{h0}& +
\ldots \cr
 \ \vdots & & &    &                     &         &           \cr
 A_r^{\bf m} & = &{A^{\bf m}_{rh} \over (s-s_0)^h} & +\ldots + & {A^{\bf m}_{r1} \over s-s_0} & +
A^{\bf m}_{r0}& + \ldots \cr  
}$$
avec $h = \min (r, o_{\bf m}(s_0))$. 

\noindent Il n'est pas difficile, en utilisant le th\'eor\`eme IV.1, de voir que toute
combinaison lin\'eaire non nulle $\sum_{i=l}^r a_iA^{\bf m}_{il}$ ($l=0, \ldots , l=r$) est \`a
support de rang $r-l$, en particulier la famille  $A^{\bf m}_{ll}, A^{\bf m}_{l+1l}, \ldots ,
A^{\bf m}_{rl}$ est lin\'eairement ind\'ependante. On montre que la distribution $T$ est
combinaison lin\'eaire des distributions $A^{\bf m}_{00},A^{\bf m}_{11}, \ldots, A^{\bf m}_{hh},
\ldots ,A^{\bf m}_{rh}$: par le lemme V.1 on a sur $V^{\times}$:
$$T = \sum a_iA^{\bf m}_{i0}.$$
La distribution $U = T- \sum a_iA^{\bf m}_{i0}$, est une distribution de
support \`a rang au plus $r-1$ et quasi-homog\`ene de degr\'e $s_0$; par le lemme V.1 sa
restriction \`a l'ouvert $W_{r-1} = \cup_0^{r-1} W_{r-1q}$ est homog\`ene, en traduisant tout
cela on montre que $a_1 = a_2 = \ldots = a_r = 0$, ce qui entraine que $U = T - a_0 A^{\bf m}_{00}$
et donc homog\`ene; nous sommes ramen\'es \`a une distribution homog\`ene de support \`a
rang au plus $r-1$; on continue ainsi de proche en proche: supposons que $h <r$, ainsi $h = o_{\bf
m}(s_0)$, et donc on a l'\'egalit\'e:
$$m_{o_{{\bf m}(s_0)}+1} -1 - (o_{{\bf m}(s_0)})d/2 < s_0 \leq m_{o_{{\bf m}(s_0)}} -1 -
(o_{{\bf m}(s_0)}-1)d/2 $$
\`a la derni\`ere \'etape, nous serons en pr\'esence d'une distribution quasi-homog\`ene \`a
support de rang au plus $r-h-1$ de la forme
$U-\sum a_iA^{\bf m}_{ih}$, mais le lemme V.1 impliquerait alors que 
$$s_0 \leq m_{o_{{\bf m}(s_0)}+1} -1 - (o_{{\bf m}(s_0)})d/2 $$
ce qui n'est pas, et par cons\'equent $U-\sum a_iA^{\bf m}_{ih} = 0$; si maintenant $h = r$, nous
arrivons \`a la derni\`ere \'etape \`a une distribution $s_0$ homog\`ene port\'ee par l'origine,
on sait qu'une telle distribution  est unique \`a un scalaire pr\`es et que donc elle est
proportionnelle \`a $A^{\bf m}_{rr}$.
Cela donne le r\'esultat pour $d$ pair.

La d\'emonstration pour le cas $d$ impair est analogue, on remplace les distributions $A_i^{\bf
m}(.,s)$ par les distributions $A_i^{\bf m}(.,s)$ ($0\leq i
\leq  E({r\over 2})$; $E(x)$ d\'esignant la partie enti\`ere de $x$) et $B_i^{\bf m}(.,s)$ ($0\leq
i
\leq  E({r-1\over 2})$) donn\'ees par:

$$A_l^{\bf m}(.,s) = \sum_{j \equiv 0 [2]}\exp (-i\pi j s_{0})j^{l}T_j^{\bf m}(,s)$$ et $$B_l^{\bf
m}(.,s) = \sum_{j
\equiv 1 [2]}\exp (-i\pi j s_{0})j^{l}T_j^{\bf m}(,s).$$

\bigskip

\bigskip

D\'esignons par ${\cal E}_{\bf m}$
l'espace vectoriel engendr\'e par les distributions $A^{\bf m}_h(.,s)$ (o\`u $s$ prend toutes les
valeurs critiques pour $\pi_{\bf m}$, $1\leq h \leq
\alpha$, et o\`u $A^{\bf
m}(.,s)$ est une combinaison lin\'eaire quelconque des $\Phi^{\bf m}_j(.,s)$) du paragraphe IV.

\proclaim Th\'eor\`eme V.4. Soit ${\bf m } = (m_1, \ldots ,m_{r-1},0)$, alors le sous espace
vectoriel des \'el\'ements \`a support singulier de ${\cal I}_{\pi_{\bf m}}$ est l'espace
${\cal E}_{\bf m}$.

Soit $T$ un
\'el\'ement non nul de ${\cal I}_{\pi_{\bf m}}$, \`a support de rang $0< p < r$ contenant
$S_{pq}$; on reprend les consid\'erations du lemme V.1, en prenant des \'el\'ements
$a$ de la forme $a = \P(o_{p,0} +
\sum_{i = p+1}^ra_ie_i)$ ($a_i$ positifs non nuls tels que $a_{p+1} \ldots a_r = 1$). Au lieu des
relations (I), on trouve alors:

$$\matrix{
n_{r-1}-n_r & = & m_{r-p-1}-m_{r-p} \cr
n_{r-2}-n_r & = & m_{r-p-2}-m_{r-p} \cr
\ \vdots  &   &                 \cr
n_{p+1}-n_r & = & m_1-m_{r-p}. \cr
} $$
 Les entiers $n_{r-1}-n_r,
\ldots , n_{p+1}-n_r$ sont donc parfaitement d\'etermin\'es par la donn\'ee de $p$ et de ${\bf m
}$, de plus l'ordre de transversalit\'e
$M$ et l'entier $n_r$ sont reli\'es par l'\'equation:

$$M -(r-p)n_r = m_1 + \ldots + m_{r-p} - (r-p)m_{r-p}.$$
Tout ceci montre comme dans le lemme V.1, que la famille
$(\psi_I(o_{p,q}))_{|I|=M}$ est unique
\`a un facteur multiplicatif  pr\`es. Prenons $s$ tel que:
$s = -(n/r-dp/2) + m_{r-p} - n_r$, on sait qu'il existe une distribution $U$ homog\`ene de degr\'e
$s$ et de support de rang $p$ contenant $S_{pq}$, le r\'esultat pr\'ec\'edent  nous montre qu'il
existe un scalaire 
$a$ telle que la distribution $T-aU$ ne
pr\'esente plus dans
$S_{pq}\cap W$ de termes d'ordre $M$, ce qui entraine que l'ordre de transversalit\'e de $T - aU$
 est inf\'erieur \`a $M$ dans $W_{pq}$ tout entier. On conclut alors par r\'ecurrence.

\bigskip

Remarque V.5: Dans [Kol.-Va.1], les auteurs font agir les op\'erateurs diff\'eren-
 \noindent tiels multiplication par $\det x$ et $\det (\partial)$ ( ces deux
op\'erateurs engendrant une alg\`ebre de Lie isomorphe ${\goth s}{\goth l}(2,
{\Bbb C})$)  sur l'espace ${\cal E}_{\bf m}$, obtenant ainsi des
${\goth s}{\goth l}(2, {\Bbb C})$-modules interessants: il apparait en particulier des
modules qui ne sont pas
\`a poids. J.A.C.Kolk et V.S. Varadarajan \'etudient
et caract\'erisent (alg\'ebriquement) ces modules. Dans le cas d'une alg\`ebre de Jordan simple et
euclidienne g\'en\'erale, l'alg\`ebre de Lie engendr\'ee par les op\'erateurs diff\'erentiels $\det
x$ et
$\det (\partial)$ n'est plus
 de dimension finie (cf.[Ru.]), mais on peut constater pour les modules ${\cal E}_{\bf m}$ des
propri\'et\'es similaires au cas \'etudi\'e dans [Kol.-Va.1]. Il serait peut-\^etre interessant
d'essayer de caract\'eriser ces modules.

\bigskip
\bigskip
\bigskip

\centerline {\bf BIBLIOGRAPHIE}

\bigskip

\noindent {\bf [Bl.1] B. Blind:} Distributions homog\`enes sur une alg\`ebre de Jordan

\noindent{\it (Bull. Soc. math. France, 125, 1997, 493-528).}

\noindent {\bf [Bl.2] B. Blind:} Fonctions zeta \`a plusieurs variables associ\'ees aux alg\`ebres
de Jordan simples euclidiennes.

\noindent{\it (C.R.Acad.Sc.Paris, 311,1990, 215-217).}

\noindent {\bf [Br.-Koe.] H.Braun, M.Koecher:} Jordan-Algebren.

\noindent{\it (Springer Verlag 1966).}

\noindent {\bf [Fa.-Ko.] J.Faraut, A.Koranyi:} Analysis on Symmetric Cones.

\noindent{\it (Clarendon Press, 1994).}

\noindent {\bf [Hi.-Ne.] J.Hilgert, K.-H.Neeb:} Vector valued Riesz distributions on Euclidean
Jordan Algebras.

\noindent{\it (Jour. Geom. Anal. 11, 2001, 43-75).}

\noindent {\bf [Ka.] S.Kaneyuki:} On the causal structures of the Shilov boundaries of symmetric
bounded domains.

\noindent{\it ( Prospects in complex geometry, Proceedings, Katata/Kyoto, 1989, J.Noguchi, T.Ohsawa
(Eds), Lec. notes in Math., 1468, Springer, 127-159 ).}

\noindent {\bf [Kol.-Va.1] J.A.C.Kolk, V.S. Varadarajan:} Lorentz invariant distributions supported
on the forward light cone.

\noindent{\it (Compositio Math. 81,1992, 61-106).}

\noindent {\bf [Kol.-Va.2] J.A.C.Kolk, V.S. Varadarajan:} Riesz Distributions

\noindent{\it (Math. Scand. 68, 1991,273-291).}

\noindent {\bf [Mu.] M.Muro:} Invariant hyperfunctions on regular prehomogeneous vector spaces of
commutative parabolic type.

\noindent{\it (T\^ohoku Math. Jour. 42, 1990, 163-193).}

\noindent {\bf [Ra.-Sc.] S. Rallis, G. Schiffmann:} Distributions invariantes par le groupe
orthogonal

\noindent{\it (Analyse Harmonique sur les Groupes de Lie S\'eminaire Nancy-Strasbourg 1973/75,
P.Eymard, J. Faraut, G. Schiffmann, R.Takahashi (Eds), Lec. notes in Math., 497,  Springer, 1975,
494-642 ).}

\noindent {\bf [Ri.-St.] F.Ricci, E.Stein:} Homogeneous distributions on spaces of Hermitean
matrices

\noindent{\it (Jour. f\"ur die reine und ang. Math. 368, 1986, 142-164).}

\noindent {\bf [Ru.] H. Rubenthaler:} Une dualit\'e de type Howe en dimension infinie

\noindent{\it (C.R.Acad.Sc.Paris, 314,1992, 435-440).}

\noindent {\bf [Sa.1] I.Satake:} Algebraic Structures of Symmetric domains.

\noindent{\it (Iwanami-Shoten and Princeton Univ. Press, 1980).}

\noindent {\bf [Sa.2] I.Satake:} On zeta functions associated with self dual homogeneous cones.

\noindent{\it (Reports on Symposium of Geometry and Automorphic Functions, T\^ohoku Univ., Sendai,
1988, 145-168).}

\noindent {\bf [Sa.-Fa.] I.Satake, J.Faraut:} The functional equation of zeta distributions
associated with formally real Jordan algebras.

\noindent{\it (T\^ohoku Math. Jour. 36, 1984, 469-482).}

\noindent {\bf [Sa.-Sh] M.Sato, T.Shintani:} On zeta function associated with prehomogeneous vector
spaces.

\noindent{\it (Ann. of Math. 100, 1974, 131-170).}

\noindent {\bf [Sp.] T.A. Springer:} Jordan Algebras and Algebraic Groups.

\noindent{\it (Springer, Ergebn. der Math. und ihrer Grenz., 75, 1973).}

\noindent {\bf [Ze.] D.P. Zelobenko:} Compact Lie Groups and their Representations.

\noindent{\it (AMS, Trans. Math. Mono. Vol 40, 1973).}

\bigskip

\bye

%% file: Init.Som.tex
\magnification =1200\hoffset=6mm\voffset=10mm


\normalbaselines
\topskip=18pt
\hsize=120mm \vsize=180mm \parindent=5mm \vskip 15mm
\font\set=cmr8\font\sept=cmti8

\font\chiffre=cmbx10 scaled 2000
\font\titre=cmbx10 scaled 1600
\def \d {\,{\rm d}}

\def \P {{\cal P}}

\def \L {{\cal L}}

\def \D {{\cal D}}

\def\epsilon{\varepsilon}

\def\\D{{\bf D}}


\catcode`\@=11

\font\tenmsx=msam10
\font\sevenmsx=msam7
\font\fivemsx=msam5
\font\tenmsy=msbm10
\font\sevenmsy=msbm7
\font\fivemsy=msbm5
\newfam\msxfam
\newfam\msyfam
\textfont\msxfam=\tenmsx  \scriptfont\msxfam=\sevenmsx
  \scriptscriptfont\msxfam=\fivemsx
\textfont\msyfam=\tenmsy  \scriptfont\msyfam=\sevenmsy
  \scriptscriptfont\msyfam=\fivemsy

\def\hexnumber@#1{\ifnum#1<10 \number#1\else
 \ifnum#1=10 A\else\ifnum#1=11 B\else\ifnum#1=12 C\else
 \ifnum#1=13 D\else\ifnum#1=14 E\else\ifnum#1=15 F\fi\fi\fi\fi\fi\fi\fi}

\def\msx@{\hexnumber@\msxfam}
\def\msy@{\hexnumber@\msyfam}
\mathchardef\nmid="3\msy@2D
\mathchardef\varnothing="0\msy@3F
\mathchardef\nexists="0\msy@40
\mathchardef\smallsetminus="2\msy@72
\def\Bbb{\ifmmode\let\next\Bbb@\else
 \def\next{\errmessage{Use \string\Bbb\space only in math mode}}\fi\next}
\def\Bbb@#1{{\Bbb@@{#1}}}
\def\Bbb@@#1{\fam\msyfam#1}

\catcode`\@=\active


\newcount\contsflag
\contsflag=0
\newcount\chstart
\def\chapcom#1#2{\vfill\eject\chstart=\pageno
\def\bookhead{ #1 \sept{#2}}%
\headline={\ifodd\pageno\rightheadline \else\leftheadline\fi}
\def \rightheadline{\ifnum \pageno=\chstart{\hfill}
\else \set {#1 \sept{#2}}\hfil\tenrm\folio \ \fi}
\def\leftheadline{\ifnum \pageno=\chstart{\hfill}
\else \tenrm \folio \  \hfil%
\set{\sept{Sets of multiples}}
\fi}}
\def\newchapter#1#2{\chapcom{#1}{#2}
\centerline {\chiffre #1}
\vskip 8mm
\centerline {\titre #2}
\vskip 20mm
\ifnum\contsflag=1{ \write0{+Contchap[#1][#2][\the\pageno]} }\fi
}
\def\newchapterl#1#2#3{\chapcom{#1}{#2}
\centerline {\chiffre #1}
\vskip 8mm
\centerline {\titre #2}
\bigskip
\centerline {\titre #3}
\vskip 20mm
\ifnum\contsflag=1{ \write0{+Contchapl[#1][#2][#3][\the\pageno]} }\fi
}
\def\Biblio{\vfill\eject\chstart=\pageno
\headline={\ifodd\pageno\rightheadline \else\leftheadline\fi}
\def \rightheadline{\ifnum \pageno=\chstart{\hfill}
\else \set {\sept Bibliography}\hfil\tenrm\folio \ \fi}
\def\leftheadline{\ifnum \pageno=\chstart{\hfill}
\else \tenrm \folio \  \hfil%
\set{Bibliography}\fi}
\centerline {\chiffre Bibliography}
\vskip 20mm
\ifnum\contsflag=1{ \write0{+Contbib[\the\pageno]} }\fi
}
%
%
\def\section#1#2{\goodbreak\noindent {\bf\S\ #1 #2}\smallskip%
\ifnum\contsflag=1{ \write0{+Contsect[#1][#2][\the\pageno]} }\fi
}
%
\def\widesection#1#2{\goodbreak\parindent=12mm
\item {\bf \S\  #1\ } {\bf #2}\smallskip
\parindent=5mm%
\ifnum\contsflag=1{ \write0{+Contsectl[#1][#2][\the\pageno]} }\fi
}
\def\Notes{\goodbreak\vskip 1in\goodbreak\centerline{\titre Notes}\vskip 0.5in%
\ifnum\contsflag=1{\write0{+Contnotes[\the\pageno]}}\fi
}
\def\Exercises{\goodbreak\vskip 1in\goodbreak\centerline{\titre Exercises}\vskip 0.5in%
\ifnum\contsflag=1{\write0{+Contexers[\the\pageno]}}\fi
}
%
%

%
%

%
%

%
%

\def\pmb#1{\setbox0=\hbox{#1}%
\kern-.025em\copy0\kern-\wd0
\kern.05em\copy0\kern-\wd0
\kern-.025em\raise .0433em\box0 }

\mathchardef\alpha="710B

%
%
 